\documentclass[preprint,12pt]{elsarticle}

\usepackage{amssymb}
\usepackage{amsmath}
\usepackage{amsfonts}
\usepackage{graphicx}
\usepackage{textcomp}
\usepackage{booktabs}
\usepackage{multirow}
\usepackage{xcolor}
\usepackage{hyperref}
\usepackage{algorithm}
\newcommand{\StateX}{\Statex\hspace{\algorithmicindent}\hspace{\algorithmicindent}}
\usepackage{algpseudocode}
\usepackage{bm}

\usepackage[margin=3cm]{geometry} 

\usepackage{tikz}

\usepackage{lineno}

\newcommand{\bft}[1]{\ooalign{#1\cr\hfil\kern.05em#1\hfil\cr\hfil\kern.02em#1\hfil}}

\journal{Electric Power Systems Research}

\begin{document}

\begin{frontmatter}




\title{End-to-End Pseudo-Measurement Learning for State Estimation under Limited Observability}



\author[uma]{J. G. De La Varga}
\author[uma]{S. Pineda}
\author[ucl]{A. Stratigakos}
\author[uma]{and J. M. Morales}

\affiliation[uma]{organization={OASYS Research Group, University of M\'alaga},
            city={M\'alaga},
            country={Spain}}
            
\affiliation[ucl]{organization={UCL Energy Institute, University College London},
            city={London},
            country={U.K.}}

\begin{abstract}
Distribution System State Estimation (DSSE) is becoming increasingly important with the integration of Distributed Energy Resources (DERs) and the active operation of distribution networks (DNs), but it remains challenging due to the limited and heterogeneous monitoring infrastructure available in these networks. 
To address this challenge, this paper proposes a novel DSSE framework that restores observability through data-driven pseudo-measurements generated by a Neural Network (NN), while preserving the exact non-linear AC power flow model within a classical Weighted Least Squares (WLS) estimator. 
Unlike conventional approaches that generate pseudo-measurements independently of the physical estimation process, the proposed method explicitly couples both components through an end-to-end learning formulation. 
Specifically, the WLS estimator is embedded as a layer within the NN architecture, enabling implicit differentiation to propagate the impact of pseudo-measurements on the final estimation error back to the NN parameters.
Extensive numerical experiments on the IEEE 30-bus and IEEE 33-bus systems demonstrate that the proposed framework consistently outperforms state-of-the-art methods in state estimation (SE) accuracy under a wide range of loading conditions and measurement configurations.
\end{abstract}




\end{frontmatter}

\section{Introduction} \label{sec:intro}

\subsection{Context, Aim and Motivation}

DSSE faces severe observability challenges due to the sparse, heterogeneous nature of current measurement infrastructures \cite{yue2024graph}. 
Many DN remain poorly monitored, hindering operators from maintaining reliability as grids transition to more active operational roles \cite{primadianto2016review}.
This data scarcity stems from a stark mismatch in sensor capabilities, as high-accuracy, fast-reporting sensors, such as SCADA systems and micro-Phasor Measurement Units ($\mu$PMUs), which offer real-time insights, are deployed in extremely limited numbers. 
Conversely, Advanced Metering Infrastructure (AMI) provides abundant smart meter data, but its high latency and slow reporting rates make it largely unsuitable for real-time tracking \cite{cheng2023survey}. 
In this context, this paper addresses the SE problem in DNs that remain unobservable when relying solely on the available physical measurements.

\subsection{Contributions}
Our work contributes to restoring network observability through the explicit generation and usage of pseudo-measurements.
We propose a SE framework built on the conventional WLS formulation, where observability is recovered through pseudo-measurements generated by a NN. 
By leveraging these data-driven pseudo-measurements, the proposed approach enables the WLS estimator to operate directly on the exact non-linear AC power flow model, thereby avoiding the structural simplifications and linear approximations commonly required in alternative methodologies.
Unlike conventional data-driven approaches, which generate pseudo-measurements \emph{independently} of the downstream physical estimation process \cite{Wang23NNpsm, Dehghanpour19game-theoretic}, 
the proposed pseudo-measurement generation mechanism is directly optimized for its final objective, namely, accurate estimation of the true network state (voltage magnitudes and angles) and the derived electrical quantities (power injections and branch power flows). 
To this end, we adopt an \emph{end-to-end} learning paradigm \cite{mandi2024decision}, in which the WLS-SE routine is embedded within the NN architecture during training as an \textit{implicit layer} \cite{agrawal2019differentiable}. 
Through implicit differentiation, this embedded optimization layer calculates the sensitivities of the final SE relative to each pseudo-measurement, allowing the NN parameters to be updated accordingly during backpropagation.
Extensive numerical experiments conducted across a wide variety of configurations (including different measurement configurations and demand variabilities) corroborate our findings and demonstrate the effectiveness and robustness of the proposed approach.



\subsection{Related Work}
Historical data on the network's past states can be used to circumvent the limitations of real-time physical telemetry.
Two main approaches exist to leverage these data for SE. 
The first approach relies on \emph{state forecasting}, where predictive models of varying complexity---ranging from Bayesian frameworks \cite{mestav2019bayesian, dobbe_linear_2020, zhang_robust_2024} to advanced artificial NN and machine learning techniques \cite{mestav2019bayesian, zhang2019real, pagnier2021physics, mukherjee2022power}---are trained on historical data to directly predict the current state of the grid from existing measurements. 
The second approach utilizes historical data to generate statistical \textit{pseudo-measurements}, 
which are then integrated into traditional SE solvers (such as WLS) to act as surrogates for the missing physical measurements \cite{dehghanpour2018survey, singh2010distribution, manitsas2012distribution, zhao2016robust, wang2020physics, zhang2023robust, vanin2023exact}.

Beyond these two predominant data-driven approaches, recent work has explored alternative mathematical and algorithmic paradigms to address the observability deficit. 
For instance, \cite{Donti2020} proposes an alternative state estimator that is based on matrix completion, specifically tailored for DNs with low observability. 
Operating as a static estimator based on a linearized power-flow model, this approach exploits the empirically observed low-rank structure of the matrix collecting system variables to reconstruct missing information from a single time instant. 
Other works focus on enhancing the underlying power-flow formulations or introducing novel optimization frameworks. 
In \cite{liu_data_2025}, a data-enhanced linear power-flow formulation is introduced, where selected entries of the admittance matrix are adjusted using available measurements. 
This method achieves robustness against model-parameter uncertainty through a Huber estimator and leverages zero-injection nodes to improve observability. 
Similarly, \cite{zeraati_novel_2024} casts the DSSE problem as a mixed-integer linear program (MILP) where both system states and nodal consumptions serve as decision variables. 
Rather than acting as direct inputs, pseudo-measurements derived from historical data are used to define feasible bounds for load variables, allowing the MILP to maximize estimation accuracy subject to power-flow constraints. 
Furthermore, efforts have been made to adapt traditional solvers to operate under extreme sparsity without relying on pseudo-measurements. 
The authors of \cite{pau_wls-based_2024} design a modified WLS estimator that enhances observability by incorporating allocation factors. 
This approach assumes that the power consumed or generated within each cluster of buses remains at the same per-unit level relative to their rated capacity, thereby retaining the advantages of a mature WLS framework even in typically unobservable scenarios. 
These alternative approaches, however, either rely on approximate models \cite{Donti2020, liu_data_2025,zeraati_novel_2024} or introduce structural simplifications that reduce the complexity of the SE~\cite{pau_wls-based_2024}.

To the best of our knowledge, the method most closely related to ours is presented in~\cite{ding2025power}, which also adopts an end-to-end learning framework by incorporating a state estimator as a differentiable layer. However, the objective and functional role of the embedded optimization layer differ substantially from those in our approach. In~\cite{ding2025power}, the differentiable estimator---built upon a convex relaxation of the AC power flow equations---is primarily used to learn adaptive measurement weights within the estimation process itself, assuming full real-time observability. 
In contrast, our methodology leverages the implicit layer to guide the learning of delayed pseudo-measurements under limited observability conditions, enforcing physical consistency directly during the forecasting stage. Therefore, rather than tuning the estimator weights, our approach integrates the network physics into the pseudo-measurement generation process, enabling improved reconstruction of both state and power variables even in scenarios with sparse sensors.

\subsection{Paper Organization}

The rest of the paper is organized as follows. Section~\ref{sec:formulation} introduces the notation, formulates the DSSE problem, presents the classical WLS estimator, and discusses observability challenges in poorly monitored networks, along with existing data-driven approaches for unobservable systems. Section~\ref{sec:il} details the proposed methodology, including pseudo-measurement generation via NNs and the integration of WLS as a differentiable implicit layer, along with the corresponding training and inference procedures. Section~\ref{sec:case-studies} reports numerical results on IEEE 30-bus and IEEE 33-bus systems under various scenarios. Finally, Section~\ref{sec:concl} concludes and outlines future work.

\section{State Estimation under Limited Observability}\label{sec:formulation}
In this section, we first introduce the notation, present the DSSE problem, the classical WLS estimator, and the observability challenges in poorly monitored networks (Subsection~\ref{sec:DSSE}).
Then, we review state forecasting (Subsection~\ref{sec:sf}) and pseudo-measurement generation (Subsection~\ref{sec:psm}).

\paragraph*{Notation} 
We use bold font lowercase (uppercase) for vectors (matrices), calligraphic for sets, $|\cdot|$ to denote set
cardinality, and $\hat{(\cdot)}$ to mark quantities produced by an estimator.
We denote the measurement vector $\mathbf{z}$ and the system state vector $\mathbf{x} = [x_1,\dots,x_n]^{\top} = [V_1,\dots,V_N,\theta_2,\dots,\theta_N]^{\top}$, where $N$ is the number of buses of the network and $n=2N-1$ (since we fix $\theta_1=0$). 
Subscript $t$ indexes time instances; superscripts $a, d$ tag \emph{available} and \emph{delayed} quantities (introduced in~\eqref{eq:za_zd_x_relation}).

\subsection{The State Estimation Problem}\label{sec:DSSE}

The relationship between the measurement vector $\mathbf{z} \in \mathbb{R}^m$ and the system state $\mathbf{x} \in \mathbb{R}^{n}$ is given by
\begin{equation}\label{eq:z_x_relation}
\mathbf{z}=\mathbf{h}(\mathbf{x})+\mathbf{e},
\end{equation}
where $\mathbf{h}(\cdot):\mathbb{R}^{n} \rightarrow \mathbb{R}^{m}$ is the measurement function which represents the power flow equations (detailed in \ref{appendix:h}),
and $\mathbf{e}\in \mathbb{R}^{m}$ is the measurement error,
which is assumed to follow a zero-mean multivariate normal distribution with independent components \cite{abur2004}.
The SE problem is formulated as a WLS optimization problem
that seeks to minimize the squared Euclidean norm of the weighted measurement errors:
\begin{equation}
\label{eq:wls-se}
\llap{\normalfont(WLS-SE)\qquad\qquad}
\hat{\mathbf{x}} \in \arg\min_\mathbf{x} \enskip
(\mathbf{z} - \mathbf{h}(\mathbf{x}))^{\top}\mathbf{W}(\mathbf{z} - \mathbf{h}(\mathbf{x})),
\end{equation}
where $\hat{\mathbf{x}}$ is the estimated state vector and $\mathbf{W}=\text{diag}\left(\bm{\sigma}^{-2}\right)$ is the weight matrix, where $\bm{\sigma}^2\in\mathbb{R}^m$ is the vector of measurement error variances, provided by the manufacturer’s noise specifications for physical measurements.  
For ease of notation, the solution of \eqref{eq:wls-se} is denoted as $\hat{\mathbf{x}}\in\texttt{wls}(\mathbf{z})$. 
In the general case, \eqref{eq:wls-se} is a non-convex, unconstrained optimization problem, 
which can be solved using an iterative Gauss--Newton or Newton--Raphson algorithm \cite{wang2019} or interior-point methods \cite{caro2020state}.  

As discussed in Section \ref{sec:intro}, even when the number and types of measurements captured by the sensors are sufficient to make the system numerically observable, 
that is, the measurement Jacobian matrix is full rank and well-conditioned~\cite{Baldwin1993ObservabilityPMU}, 
some measurements may be unavailable in real time due to delays or communication issues. 
Following \cite{yue2024graph, delaVarga2025}, we categorize distribution system measurements into two groups: 
\emph{(i)} fast and sparse measurements obtained from PMUs and SCADA, and \emph{(ii)} slow but abundant measurements obtained from Smart Meters.
Fast and sparse measurements are considered \emph{available} during real-time operations, 
whereas slow but abundant ones are typically \emph{delayed} and, hence, cannot be used for SE. 
Consequently, we partition the measurement index set as $\{1,\ldots,m\} = \mathcal{M}^a \cup \mathcal{M}^d$ into available  and delayed measurements, with $m_a := |\mathcal{M}^a|$ and $m_d := |\mathcal{M}^d|$. We denote by $\mathbf{z}^a \in \mathbb{R}^{m_a}$ and $\mathbf{z}^d \in \mathbb{R}^{m_d}$ the corresponding subvectors of $\mathbf{z}$ obtained by selecting the components indexed by $\mathcal{M}^a$ and $\mathcal{M}^d$, respectively.
The measurement function splits accordingly:
\begin{subequations}\label{eq:za_zd_x_relation}
\begin{align}
& \mathbf{z}^a=\mathbf{h}^a(\mathbf{x})+\mathbf{e}^a, \\
& \mathbf{z}^d=\mathbf{h}^d(\mathbf{x})+\mathbf{e}^d,
\end{align}
\end{subequations}
where $\mathbf{h}^a(\cdot):\mathbb{R}^{n} \rightarrow \mathbb{R}^{m_a}$, $\mathbf{h}^d(\cdot):\mathbb{R}^{n} \rightarrow \mathbb{R}^{m_d}$,
$\mathbf{e}^a\in \mathbb{R}^{m_a}$, and $\mathbf{e}^d\in \mathbb{R}^{m_d}$.
In this work, we assume that the number of measurements available in real time is insufficient to render the system observable \cite{cheng2023survey, gomez-exposito2015}, i.e., $m_a < n$. Notwithstanding this limitation, a retrospective SE problem can always be solved by assuming that all measurements become available after the fact, that is, we assume the set of available and delayed measurements makes the system observable. This retrospective system state is denoted as $\hat{\mathbf{x}}^d$ and can be computed via WLS as $\hat{\mathbf{x}}^d \in \texttt{wls}(\mathbf{z})$, where $\mathbf{z} := \big[\mathbf{z}^a, \mathbf{z}^d\big]$.

Building upon this retrospective formulation, a historical dataset $\mathcal{D} = \{(\mathbf{z}^a_t, \mathbf{z}^d_t, \hat{\mathbf{x}}^d_t)\}_{t \in \mathcal{T}}$ can be constructed. Because operators only have access to a limited set of physical measurements during real-time operation, data-driven methodologies have emerged to utilize these historical records to mitigate the lack of online observability. 
In the following, we review the two predominant data-driven paradigms in the DSSE literature, 
namely \emph{state forecasting} and \emph{pseudo-measurement generation}~\cite{delaVarga2025}.

\subsection{State Forecasting (SF)}\label{sec:sf}
 
The SF strategy bypasses the physical model entirely and learns a direct map from the available measurements to the system state. 
Letting $f_{\bm{\phi}}: \mathbb{R}^{m_a} \to \mathbb{R}^{n}$ denote the regression model parameterized by $\bm{\phi}$, the training problem reads
\begin{equation}\label{eq:sf_loss}
\hat{\bm{\phi}}^{\mathrm{SF}} \in \arg\min_{\bm{\phi}} \frac{1}{|\mathcal{T}|} \sum_{t \in \mathcal{T}} \big\| f_{\bm{\phi}}(\mathbf{z}^a_t) - \hat{\mathbf{x}}^d_t \big\|_2^2,
\end{equation}
where the supervisory signal is the retrospective WLS estimate $\hat{\mathbf{x}}^d_t$ present in $\mathcal{D}$. At inference, the predicted state at instance $t'$ is simply
\begin{equation}\label{eq:sf_inference}
\hat{\mathbf{x}}^{\mathrm{SF}}_{t'} = f_{\hat{\bm{\phi}}^{\mathrm{SF}}}(\mathbf{z}^a_{t'}).
\end{equation}
The SF strategy has the simplest deployment pipeline, as it only requires a single forward pass through $f_{\hat{\bm{\phi}}^{\mathrm{SF}}}$ to produce the state estimate. 
Its main limitation is that the estimated state need not satisfy any of the physical relations encoded in $h$: the residual $\mathbf{z}_t - \mathbf{h}(\hat{\mathbf{x}}^{\mathrm{SF}}_t)$ is not minimized in any sense. 
Equivalently, SF treats the entire state-estimation pipeline as a black-box regression problem, discarding the problem's physical structure.


\subsection{Pseudo-Measurement Generation (PS)}\label{sec:psm}
 
PS, which follows the dominant paradigm in the DSSE, learns a surrogate for the delayed measurements and feeds the resulting full measurement vector into the WLS estimator \eqref{eq:wls-se}.
We consider $g_{\bm{\phi}}: \mathbb{R}^{m_a} \to \mathbb{R}^{m_d}$ a feed-forward NN whose training problem is the supervised reconstruction error of the delayed measurements,
\begin{equation}\label{eq:psm_loss}
\hat{\bm{\phi}}^{\mathrm{PS}} \in \arg\min_{\bm{\phi}} \frac{1}{|\mathcal{T}|} \sum_{t \in \mathcal{T}} \big\| g_{\bm{\phi}}(\mathbf{z}^a_t) - \mathbf{z}^d_t \big\|_2^2.
\end{equation}
At inference, the predicted pseudo-measurements $\hat{\mathbf{z}}^d_{t'} = g_{\hat{\bm{\phi}}^{\mathrm{PS}}}(\mathbf{z}^a_{t'})$ are concatenated with the available real-time measurements $\mathbf{z}^a_{t'}$, and the resulting vector is fed into \eqref{eq:wls-se},
\begin{equation}\label{eq:psm_inference}
\hat{\mathbf{x}}^{\mathrm{PS}}_{t'} \in \texttt{wls}\big(\big[\mathbf{z}^a_{t'},\, \hat{\mathbf{z}}^d_{t'}\big]\big).
\end{equation}
 
A natural choice for the WLS weight matrix $\mathbf{W}$ in~\eqref{eq:psm_inference} combines the manufacturer-provided sensor noise specifications for physical measurements with an empirical estimate of the pseudo-measurement uncertainty. Concretely, for each delayed pseudo-measurement we set the corresponding diagonal entry of $\mathbf{W}$ to $\left(\hat{\bm{\sigma}}^d\right)^{-2}$, with
\begin{equation}\label{eq:psm_weights}
\left(\hat{\bm{\sigma}}^d\right)^{2} = \frac{1}{|\mathcal{T}|} \sum_{t \in \mathcal{T}} \Big( g_{\hat{\bm{\phi}}^{\mathrm{PS}}}(\mathbf{z}^a_t) - \mathbf{z}^d_t \Big)^2,
\end{equation}
while the weights associated with the available measurements, $\left(\hat{\bm{\sigma}}^a\right)^{-2}$, follow directly from the same sensor specifications. 

In contrast to SF, PS preserves the physical consistency of the estimate, since $\hat{\mathbf{x}}^{\mathrm{PS}}_{t'}$ is constructed as the WLS solution under the full nonlinear AC model. However, a key limitation arises because the predictor $g_{\bm{\phi}}$ is trained \emph{agnostically} with respect to this downstream estimation task. 
Because \eqref{eq:psm_loss} treats all components of $\mathbf{z}^d$ symmetrically, it minimizes error purely in the measurement space, ignoring the fact that the state variables exhibit highly asymmetric sensitivities to different measurements. 
Consequently, while the training objective of $g_{\bm{\phi}}$ is tailored to maximize the accuracy of the pseudo-measurements themselves, the actual downstream interest lies in the accuracy of the resulting state variables, creating a fundamental misalignment between the training loss and the operational goal.

\section{End-to-end Pseudo-measurement Learning}
\label{sec:il}

\begin{figure}[tb!]
    \centering
    \resizebox{\columnwidth}{!}{
    \input{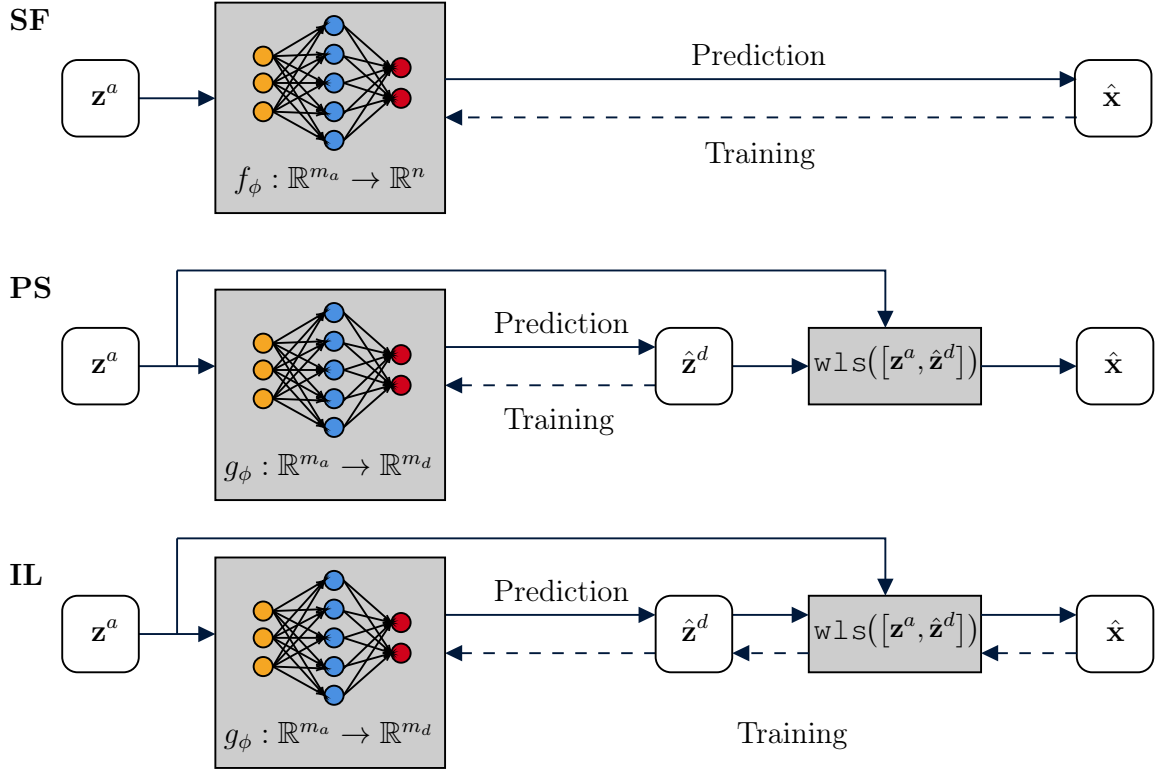}}
    \caption{Illustration of the three data-driven strategies State Forecasting (SF), Pseudo-measurements (PS), and the Implicit Layer (IL) proposed framework; contrasting their data flow at inference (solid line) with their gradient flow at training (dashed line).}
\label{fig:flowchart}
\end{figure}

SF and PS occupy opposite corners of the same trade-off: SF aligns training with the SE target but discards the physical model, while PS preserves the physics but trains its predictor without taking the state into account. The methodology proposed in this section retains the WLS solve of PS while replacing its training objective with one that propagates the SE error back to the pseudo-measurement generator. Since this end-to-end coupling is realized by embedding the WLS estimator as an implicit differentiable layer, we denote our proposal IL.
Figure~\ref{fig:flowchart} compares IL with the traditional SF and PS paradigms. 

We first formalize the IL approach as a bilevel optimization problem with a hybrid loss (in Subsection~\ref{sec:il_formulation}). 
Then, we develop the solution methodology: implicit differentiation of the WLS-SE, the end-to-end training loop, and the inference pipeline (in Subsection~\ref{sec:il_methodology}). 


\subsection{Proposed Approach Formulation}
\label{sec:il_formulation}

The IL strategy also utilizes the NN function $g_{\bm{\phi}}: \mathbb{R}^{m_a} \to \mathbb{R}^{m_d}$ as a pseudo-measurement generator, but with a different training loss. Formally, the parameters of the function $g_{\bm{\phi}}$ are the solution to the following bilevel optimization problem
\begin{subequations}\label{eq:il_bilevel}
\begin{align}
\hat{\bm{\phi}}^{\mathrm{IL}} \in \arg\min_{\bm{\phi}} \quad & \mathcal{L}(\bm{\phi},\gamma) \label{eq:il_outer} \\
\text{s.t.} \quad & \hat{\mathbf{x}}_t(\bm{\phi}) \in \texttt{wls}\big(\big[\mathbf{z}^a_t,\, \hat{\mathbf{z}}^d_t(\bm{\phi})\big]\big),\label{eq:il_inner}\\
& \hat{\mathbf{z}}^d_t(\bm{\phi}) = g_{\bm{\phi}}(\mathbf{z}^a_t), \quad \forall t \in \mathcal{T},
\end{align}
\end{subequations}
where the inner problem~\eqref{eq:il_inner} is the WLS estimator instantiated with the augmented measurement vector $\big[\mathbf{z}^a_t,\, \hat{\mathbf{z}}^d_t(\bm{\phi})\big]$, and the outer objective $\mathcal{L}(\bm{\phi},\gamma)$ is the convex combination
\begin{equation}\label{eq:il_hybrid_loss}
\mathcal{L}(\bm{\phi},\gamma) \;=\; \frac{1}{|\mathcal{T}|} \sum_{t \in \mathcal{T}} \Big[ \gamma\,\|\hat{\mathbf{x}}_t(\bm{\phi}) - \hat{\mathbf{x}}^d_t\|_2^2 + (1-\gamma)\,\|\hat{\mathbf{z}}^d_t(\bm{\phi}) - \mathbf{z}^d_t\|_2^2 \Big].
\end{equation}
The hyperparameter $\gamma \in [0,1]$ controls the trade-off between two complementary training signals. At $\gamma = 0$, IL reduces to PS and~\eqref{eq:il_bilevel} collapses to the conventional decoupled pipeline in which the pseudo-measurement predictor is trained independently of the downstream WLS. At $\gamma = 1$, training becomes fully end-to-end as $g_{\boldsymbol{\phi}}$ is optimized solely so that the resulting state estimate $\hat{\mathbf{x}}_t(\boldsymbol{\phi})$ matches the reference state $\hat{\mathbf{x}}^d_t$, with no requirement that the outputs of the NN $g_{\boldsymbol{\phi}}$ themselves, i.e., the reconstructed delayed measurements, resemble plausible measurement values. 
For $\gamma\in(0,1)$, the task term aligns the estimation of the NN parameter vector  $\boldsymbol{\phi}$ with the minimization of the state error, while the regularization term anchors $g_{\boldsymbol{\phi}}$ to physically meaningful pseudo-measurements. 
In other words, since many distinct pseudo-measurement assignments can yield comparable state estimates through the inner WLS, the regularization term plays an additional disambiguating role by selecting the predictions closest to the observed $\mathbf{z}^d_t$ amongst the ones that are consistent with a good state estimate.

\subsection{Solution Methodology}
\label{sec:il_methodology}

Solving~\eqref{eq:il_bilevel} by gradient methods requires the Jacobian of \eqref{eq:il_inner} with respect to its measurement input whenever $\gamma > 0$. We first obtain it by implicit differentiation, 
avoiding the need to backpropagate through the unrolled iterates of the WLS solver. Then, we define the end-to-end training loop of the gradient-based method. Finally, the inference pipeline of the methodology at deployment is outlined.

\subsubsection{Implicit Differentiation of WLS-SE}\label{sec:imp_diff}

For the unconstrained non-linear least-squares problem~\eqref{eq:wls-se}, we can derive the first-order optimality condition at the solution $\hat{\mathbf{x}}$, which reads
\begin{equation}\label{eq:foc}
\mathbf{F}(\hat{\mathbf{x}},\, \mathbf{z}) \;:=\; -\mathbf{J}(\hat{\mathbf{x}})^{\!\top} \mathbf{W} \big(\mathbf{z} - \mathbf{h}(\hat{\mathbf{x}})\big) \;=\; \mathbf{0},
\end{equation}
where $\mathbf{J}(\hat{\mathbf{x}}) := \partial \mathbf{h} / \partial \mathbf{x} \,\big|_{\hat{\mathbf{x}}} \in \mathbb{R}^{m\times n}$ is the measurement Jacobian evaluated at $\hat{\mathbf{x}}$. Equation~\eqref{eq:foc} implicitly defines 
a vector function $\hat{\mathbf{x}}(\mathbf{z})$ such that $\mathbf{F}(\hat{\mathbf{x}}(\mathbf{z}),\, \mathbf{z}) = 0$. 
Then, provided that the Jacobian is non-singular at the solution---a condition equivalent to the system being numerically observable when the augmented vector $\big[\mathbf{z}^a, \hat{\mathbf{z}}^d\big]$ is fed to the WLS---the implicit function theorem yields
\begin{equation}\label{eq:imp_difg_general}
\frac{\partial \hat{\mathbf{x}}}{\partial \mathbf{z}} \;=\; -\Big(\frac{\partial \mathbf{F}}{\partial \hat{\mathbf{x}}}\Big)^{\!-1} \frac{\partial \mathbf{F}}{\partial \mathbf{z}}.
\end{equation}
Differentiating~\eqref{eq:foc} explicitly, first with respect to the state, 
\begin{equation}\label{eq:dg_dx}
\frac{\partial \mathbf{F}}{\partial \hat{\mathbf{x}}} \;=\; \mathbf{J}(\hat{\mathbf{x}})^{\!\top} \mathbf{W} \mathbf{J}(\hat{\mathbf{x}}) \;-\; \sum_{k=1}^{m} \big[\mathbf{W}\,\mathbf{r}(\hat{\mathbf{x}},\mathbf{z})\big]_k \, \nabla^2_{\mathbf{x}} h_k(\hat{\mathbf{x}}),
\end{equation}
with residual $\mathbf{r}(\hat{\mathbf{x}},\mathbf{z}) = \mathbf{z} - \mathbf{h}(\hat{\mathbf{x}})$, and then with respect to the measurements,
\begin{equation}\label{eq:dg_dz}
\frac{\partial \mathbf{F}}{\partial \mathbf{z}} \;=\; -\,\mathbf{J}(\hat{\mathbf{x}})^{\!\top} \mathbf{W}.
\end{equation}
Following standard practice in WLS SE~\cite{abur2004}, we adopt the Gauss--Newton (GN) approximation, which discards the second-derivative term in~\eqref{eq:dg_dx}. Substituting the GN approximation $\partial \mathbf{F}/\partial \hat{\mathbf{x}} \approx \mathbf{J}^{\!\top} \mathbf{W} \mathbf{J}$ (where the dependence of $\mathbf{J}$ on $\hat{\mathbf{x}}$ has been left implicit) together with~\eqref{eq:dg_dz} into~\eqref{eq:imp_difg_general} yields the closed-form sensitivity
\begin{equation}\label{eq:wls_jacobian}
\frac{\partial \hat{\mathbf{x}}}{\partial \mathbf{z}} \;=\; \big(\mathbf{J}^{\!\top} \mathbf{W} \mathbf{J}\big)^{\!-1} \mathbf{J}^{\!\top} \mathbf{W}.
\end{equation}
Equation~\eqref{eq:wls_jacobian} is the well-known sensitivity matrix of the WLS estimator~\cite{abur2004} and recovers, in the linear case $\mathbf{h}(\mathbf{x}) = \mathbf{H}\mathbf{x}$, the standard linear WLS pseudo-inverse $(\mathbf{H}^{\!\top} \mathbf{W} \mathbf{H})^{-1} \mathbf{H}^{\!\top} \mathbf{W}$.

\subsubsection{End-to-End Training Loop}\label{sec:training_loop}

\begin{algorithm}[!ht]
\caption{End-to-End Training of IL}
\label{alg:il_training}
\begin{algorithmic}[1]
\Require Historical dataset $\mathcal{D} = \{(\mathbf{z}^a_t, \mathbf{z}^d_t, \hat{\mathbf{x}}^d_t)\}_{t \in \mathcal{T}}$, weight matrix $\mathbf{W}$, learning rate $\eta$, max epochs $E$, GN iterations $K$, loss hyperparameter $\gamma \in [0,1]$, early stopping patience $\Psi$ and threshold $\varphi$
\State Split $\mathcal{D} = \mathcal{D}^{\mathrm{train}} \cup \mathcal{D}^{\mathrm{val}}$
\State Initialize $\bm{\phi}^{0}$; set $\hat{\bm{\phi}} \gets \bm{\phi}^{0}$, $\hat{\mathcal{L}} \gets +\infty$, $\psi \gets 0$
\For{$e = 1, \ldots, E$}
    \For{minibatch $\mathcal{B} \subset \mathcal{D}^{\mathrm{train}}$}
        \State \textbf{Forward (NN):} 
        \StateX $\hat{\mathbf{z}}^d_t \gets g_{\bm{\phi}}(\mathbf{z}^a_t)$  \quad $\forall t \in \mathcal{B}$
        \State \textbf{Forward (WLS):} 
        \StateX $\hat{\mathbf{x}}_t \gets \texttt{wls}\big(\big[\mathbf{z}^a_t,\, \hat{\mathbf{z}}^d_t\big]\big)$ via $K$ GN iterations; cache $\mathbf{J}(\hat{\mathbf{x}}_t)$
        \State \textbf{Loss:} 
        \StateX $\mathcal{L}(\bm{\phi},\gamma) \gets \tfrac{1}{|\mathcal{B}|} \sum_{t \in \mathcal{B}} \!\big[ \gamma\,\|\hat{\mathbf{x}}_t - \hat{\mathbf{x}}^d_t\|_2^2 + (1{-}\gamma)\,\|\hat{\mathbf{z}}^d_t - \mathbf{z}^d_t\|_2^2 \big]$
        \State \textbf{Backward (WLS):}
        \StateX solve $\bigl(\mathbf{J}(\hat{\mathbf{x}}_t)^{\!\top} \mathbf{W} \mathbf{J}(\hat{\mathbf{x}}_t)\bigr)\, \bm{\lambda}_t = \nabla_{\hat{\mathbf{x}}} \mathcal{L} \quad \forall t \in \mathcal{B}$
        \StateX $\big(\partial \hat{\mathbf{x}} / \partial \mathbf{z}\big)^{\!\top} \nabla_{\hat{\mathbf{x}}} \mathcal{L} \gets \mathbf{W}\, \mathbf{J}(\hat{\mathbf{x}}_t)\, \bm{\lambda}_t$
        \State \textbf{Backward (NN):} 
        \StateX compute $\partial \mathbf{z} / \partial \bm{\phi}$ via reverse-mode automatic differentiation
        \StateX $\nabla_{\mathbf{z}} \mathcal{L} \gets 2(1-\gamma)(\hat{\mathbf{z}}^d_t - \mathbf{z}^d_t) / |\mathcal{B}| \quad \forall t \in \mathcal{B}$
        \State \textbf{Parameter Update:}
        \StateX compute $\nabla_{\bm{\phi}} \mathcal{L}$ using \eqref{eq:loss_gradient}
        \StateX $\bm{\phi} \gets \bm{\phi} - \eta \, \nabla_{\bm{\phi}}\mathcal{L}$
    \EndFor
    \State \textbf{Validation:} evaluate $\mathcal{L}_{\mathrm{val}} \gets \mathcal{L}\bigl(\bm{\phi}, \gamma;\, \mathcal{D}^{\mathrm{val}}\bigr)$
    \If{$\mathcal{L}_{\mathrm{val}} < \hat{\mathcal{L}} - \varphi$}
        \State $\hat{\bm{\phi}} \gets \bm{\phi}$, \; $\hat{\mathcal{L}} \gets \mathcal{L}_{\mathrm{val}}$, \; $\psi \gets 0$
    \Else
        \State $\psi \gets \psi + 1$
    \EndIf
    \If{$\psi \geq \Psi$} \textbf{break} \EndIf
\EndFor
\State \Return $\hat{\bm{\phi}}^{\mathrm{IL}} \gets \hat{\bm{\phi}}$
\end{algorithmic}
\end{algorithm}

Embedding the WLS-SE \eqref{eq:wls-se} as an implicit layer inside the NN amounts to executing the gradient-based training procedure summarized in Algorithm~\ref{alg:il_training}. 
Most of its structure (train/validation splitting, mini-batch sampling, forward pass through the predictor, evaluation of the task loss on the validation set, parameter update, early stopping with a patience threshold) is standard supervised learning, with two key exceptions. 
First, the forward pass embeds a WLS solve via $K$ GN iterations (a user-defined parameter) between the predictor output and the loss (lines~4--6). 
Second, the backward pass must propagate the gradient \emph{through} that WLS solution. 
The remainder of this subsection derives the latter step and shows that it admits a closed form whose cost is independent of $K$.

Adopting the convention that all gradients are column vectors, the chain rule applied to the hybrid loss~\eqref{eq:il_hybrid_loss} decomposes the gradient $\nabla_{\bm{\phi}} \mathcal{L}$ into two contributions, corresponding to the two paths along which the predictor parameters $\bm{\phi}$ influence $\mathcal{L}$:
\begin{equation} 
\label{eq:loss_gradient}
\nabla_{\bm{\phi}} \mathcal{L} \;=\; \biggl(\frac{\partial \mathbf{z}}{\partial \bm{\phi}}\biggr)^{\!\top} \biggl(\frac{\partial \hat{\mathbf{x}}}{\partial \mathbf{z}}\biggr)^{\!\top} \nabla_{\hat{\mathbf{x}}} \mathcal{L} \;+\; 
\biggl(\frac{\partial \mathbf{z}}{\partial \bm{\phi}}\biggr)^{\!\top} \nabla_{\mathbf{z}} \mathcal{L}.
\end{equation}
In the first term, the upstream gradient $\nabla_{\hat{\mathbf{x}}} \mathcal{L} \in \mathbb{R}^n$ goes through the middle factor $\partial \hat{\mathbf{x}}/\partial \mathbf{z} \in \mathbb{R}^{n \times m}$, which is the WLS sensitivity~\eqref{eq:wls_jacobian}, and finally into the predictor Jacobian $\partial \mathbf{z}/\partial \bm{\phi}$. In the second term, the upstream gradient $\nabla_{\mathbf{z}} \mathcal{L} \in \mathbb{R}^{m}$ goes directly through $g_{\bm{\phi}}$ into the predictor Jacobian $\partial \mathbf{z}/\partial \bm{\phi}$. Note that $\partial \mathbf{z}/\partial \bm{\phi}$ simplifies to $\partial \mathbf{z}^d/\partial \bm{\phi}$ since $\mathbf{z}^a$ is constant. All factors are handled by reverse-mode automatic differentiation since they are composed of elementary differentiable operations. However, instead of computing 
\begin{equation}\label{eq:vjp_full}
\biggl(\frac{\partial \hat{\mathbf{x}}}{\partial \mathbf{z}}\biggr)^{\!\top}\!\nabla_{\hat{\mathbf{x}}} \mathcal{L} \;=\; \mathbf{W}\,\mathbf{J}\,\bigl(\mathbf{J}^{\!\top} \mathbf{W} \mathbf{J}\bigr)^{-1}\,\nabla_{\hat{\mathbf{x}}} \mathcal{L},
\end{equation}
we introduce an auxiliary vector $\bm{\lambda} \in \mathbb{R}^n$, defined as the solution of the $n \times n$ symmetric positive-definite system
\begin{equation}\label{eq:adjoint_solve}
\bigl(\mathbf{J}^{\!\top} \mathbf{W} \mathbf{J}\bigr)\,\bm{\lambda} \;=\; \nabla_{\hat{\mathbf{x}}} \mathcal{L},
\end{equation}
after which the downstream gradient is recovered by a single matrix-vector product
\begin{equation}\label{eq:adjoint_z}
\biggl(\frac{\partial \hat{\mathbf{x}}}{\partial \mathbf{z}}\biggr)^{\!\top} \nabla_{\hat{\mathbf{x}}} \mathcal{L} \;=\; \mathbf{W}\,\mathbf{J}\,\bm{\lambda} \;\in\; \mathbb{R}^m.
\end{equation}
The $n$-dimensional vector $\bm{\lambda}$ is the so-called \emph{adjoint} variable of the WLS optimality condition~\cite{nocedal2006numerical}, 
which defines how the scalar loss responds to infinitesimal perturbations of the equation $\mathbf{J}^{\!\top} \mathbf{W} \bigl(\mathbf{z} - \mathbf{h}(\hat{\mathbf{x}})\bigr) = \mathbf{0}$ that defines $\hat{\mathbf{x}}$ implicitly. Line~8 of Algorithm~\ref{alg:il_training} carries out exactly~\eqref{eq:adjoint_solve}--\eqref{eq:adjoint_z}.

The cost of the entire backward pass through the WLS layer therefore reduces to solving a $n \times n$ system of linear equations
and one matrix-vector product, \emph{independent of the number of GN iterations $K$ used in the forward pass}. Thus, with implicit differentiation, only $\hat{\mathbf{x}}$ and $\mathbf{J}$ at the converged solution need to be cached.  
We use the same $\mathbf{W}$ for IL as for PS, with diagonal entries given by~\eqref{eq:psm_weights}, so that PS and IL differ \emph{only} in the training loss---not in the inference pipeline, the model size, or the weighting strategy.


 
At deployment, IL inherits the inference pipeline of PS, with identical computational cost:
\begin{equation}\label{eq:il_inference}
\hat{\mathbf{z}}^d_{t'} = g_{\hat{\bm{\phi}}^{\mathrm{IL}}}(\mathbf{z}^a_{t'}), \qquad
\hat{\mathbf{x}}^{\mathrm{IL}}_{t'} \in \texttt{wls}\big(\big[\mathbf{z}^a_{t'}, \hat{\mathbf{z}}^d_{t'}\big]\big).
\end{equation}

\section{Case Studies}\label{sec:case-studies}
In this section, we first describe the experimental setup developed to evaluate the different SE methods (Subsection~\ref{sec:exp}). 
Next, we analyze the results of the numerical experiments (Subsection~\ref{sec:results}).

\subsection{Experimental Setup}\label{sec:exp}

\subsubsection{Case studies and data}
We evaluate the proposed approach considering two DNs of comparable size but with distinct topological characteristics. 
First, we consider the IEEE 30-bus 
standard benchmark for ADNs, characterized by a meshed configuration and the integration of distributed generation \cite{case30asADN_1, case30asADN_2}. 
Second, we consider the IEEE 33-bus system as a passive and radial DN \cite{BaranWu89ieee33}, which lacks distributed generation and relies solely on substation injections. 
%
We evaluate robustness considering six scenarios across the different grids. 
For the IEEE 30-bus system, we assess the performance of the proposed end-to-end learning approach in a meshed network as real-time monitoring decreases by varying the number of available PMUs in the real-time measurement set $\mathbf{z}^a$, leading to three scenarios: HIG, MED, and LOW. 
In these cases, the delayed measurement set $\mathbf{z}^d$ comprises power injections at all buses and power flows through all distribution lines. 
In contrast, the IEEE 33-bus radial system utilizes topology-driven scenarios to reflect practical utility strategies. 
The Bifurcation (BIF) scenario strategically places sensors at major branch splits to monitor power distribution into laterals, while the End-of-Line (END) scenario targets leaf nodes to capture the most volatile voltage boundary conditions. 
Finally, the PMU scenario evaluates the impact of a single high-precision measurement at a strategic midpoint. 
For this passive network, $\mathbf{z}^d$ consists of power injections at all buses, representing consumption meters at the customer nodes. 
Details of the measurement set configurations across all scenarios are provided in Table~\ref{tab:scenarios}.


\begin{table}[h!]
\centering
\caption{Measurement Configurations}
\label{tab:scenarios}
\small
\begin{tabular}{@{}llll@{}}
\toprule
\textbf{Network} & \textbf{Scenario}
  & \textbf{Available ($\mathbf{z}^a$)}
  & \textbf{Delayed ($\mathbf{z}^d$)} \\
\midrule
\multirow{6}{*}{IEEE 30}
  & HIG
  & \shortstack[l]{$V$ at buses $\{1, 2, 5, 6, 8, 13, 19, 22, 28\}$ \\
                   $\theta$ at buses $\{2, 5, 6, 8, 13, 19, 22, 28\}$}
  & \multirow{6}{*}{\shortstack[l]{$P, Q$ injections \\ at all buses \\
                                   $P, Q$ flows \\ through all branches}} \\
\cmidrule(lr){2-3}
  & MED
  & \shortstack[l]{$V$ at buses $\{1, 2, 5, 6, 8, 28\}$ \\
                   $\theta$ at buses $\{2, 5, 6, 8, 28\}$}
  & \\
\cmidrule(lr){2-3}
  & LOW
  & \shortstack[l]{$V$ at buses $\{1, 12, 21\}$ \\
                   $\theta$ at buses $\{12, 21\}$}
  & \\
\midrule
\multirow{6}{*}{IEEE 33}
  & BIF
  & \shortstack[l]{Substation $V, P, Q$ (Bus 1, Branch 1) \\
                   $P, Q$ flows at bifurcation branches $\{18, 22, 25\}$}
  & \multirow{6}{*}{\shortstack[l]{$P, Q$ injections \\ at all buses}} \\
\cmidrule(lr){2-3}
  & END
  & \shortstack[l]{Substation $V, P, Q$ (Bus 1, Branch 1) \\
                   $V$ at leaf nodes $\{18, 22, 25, 33\}$}
  & \\
\cmidrule(lr){2-3}
  & PMU
  & \shortstack[l]{Substation $V, P, Q$ (Bus 1, Branch 1) \\
                   Mid-point PMU ($V, \theta$) at bus $\{6\}$}
  & \\
\bottomrule
\end{tabular}
\end{table}

For each network, we generate demand scenarios by perturbing nominal values using uniform distributions with variabilities of 5\%, 10\%, and 20\%. 
For the 30-bus system, we obtain operating conditions by performing an Optimal Power Flow for each demand profile; 
for the 33-bus system, a standard power flow is utilized, as the substation serves as the sole generating source. 
A dataset of 10\,000 samples is generated for each combination of network and variability, partitioned into 7\,000 samples for training, 1\,500 for validation, and 1\,500 for testing. 
The overall number of configurations is 18 (two networks, three sensor scenarios, and three levels of demand variability).
To ensure statistical robustness and mitigate the impact of data partitioning, each configuration is executed five times using different random splits for the training, validation, and test sets. 
Final performance metrics are reported as the average and standard deviation across these five runs.

To simulate realistic sensor inaccuracies, additive white Gaussian noise is applied to all measurements. 
Following standard industry profiles, voltage magnitude measurements ($V$) are modeled with a standard deviation of 0.001 p.u., while voltage angle measurements ($\theta$) from PMUs are subject to a noise standard deviation of $0.1^\circ$. 
For power measurements, both nodal injections ($P, Q$) and branch flows ($P_{f}, Q_{f}$) are modeled with a fixed noise standard deviation of 0.01 p.u.

\subsubsection{Hyperparameters}
The proposed NN architectures of $f_{\bm{\phi}} \text{ and } g_{\bm{\phi}}$ consist of a fully connected multilayer perceptron with three hidden layers, each containing 128 neurons and utilizing the Rectified Linear Unit (ReLU) activation function. To ensure robust training, we utilize the AdamW optimizer with an initial learning rate $\eta=10^{-4}$ and a weight decay of $10^{-5}$ for regularization. A batch size of 1\,000 samples is employed across a maximum of 50\,000 training epochs. To manage the learning dynamics, we implement a learning rate scheduler (\texttt{ReduceLROnPlateau}) that reduces the learning rate by a factor of 0.5 if the validation loss fails to improve for 50 consecutive epochs. 
Additionally, an early stopping mechanism is integrated with a patience of $\Psi = 100$ epochs (monitored via validation loss with a threshold of $\varphi =10^{-7}$) to prevent overfitting and ensure the model converges to the most stable physical state. 

The parameters of the implicit-layer network are warm-started from the weights of the pre-trained PS network, with Gaussian noise of standard deviation $10^{-2}$ added element-wise to each parameter tensor. 
This initialization strategy provides a stable and well-conditioned starting point for the bilevel fine-tuning stage while breaking the exact symmetry of the PS solution, helping the optimizer escape potential saddle points.

\subsubsection{Benchmarks and metrics}
We evaluate the performance of the three previously introduced methodologies: SF, PS, and the IL proposed framework. 
We evaluate the sensitivity of IL to the loss hyperparameter $\gamma$ by testing three variants---IL1, IL5, and IL9---corresponding to $\gamma \in \{0.1, 0.5, 0.9\}$, respectively. In order to assess their accuracy in estimating system state variables and derived physical quantities, we compute the Root Mean Square Error (RMSE) across six key metrics: voltage magnitudes and angles (state variables), as well as active and reactive power injections and branch power flows (derived magnitudes).

\subsection{Results}\label{sec:results}

Table~\ref{tab:summary_case30} reports the mean performance across all demand variability levels, sensor configurations, and five independent stochastic runs for the 30-bus ADN, while Table~\ref{tab:summary_case33} equivalently reports for the IEEE 33-bus passive radial system. 
For each metric, the best-performing model (lowest RMSE) is indicated in bold; when mean values are identical, the model with the lowest standard deviation is considered the best.
Detailed disaggregated results across specific variability levels and sensing scenarios are provided in Table~\ref{tab:results_case_30} and  Table~\ref{tab:results_case33} in~\ref{appendix:detailed_results}.

\subsubsection{IEEE 30-bus active DN}


\begin{table*}[!t]
\caption{Average results for IEEE 30-bus (mean $\pm$ std) for different variability levels, measurement configurations, and data training sets.}
\label{tab:summary_case30}
\centering
\scriptsize
\setlength{\tabcolsep}{2pt}
\newcommand{\na}{---}
\resizebox{\textwidth}{!}{
\renewcommand{\arraystretch}{1.2}
\begin{tabular}{ccccccc}
\toprule
 & \textbf{RMSE$_{V}$} & \textbf{RMSE$_{\theta}$} & \textbf{RMSE$_{P}$} & \textbf{RMSE$_{Q}$} & \textbf{RMSE$_{P_f}$} & \textbf{RMSE$_{Q_f}$} \\
\midrule
IL1 & 1.38e-3$\!\pm\!$3.26e-4 & 1.62e-3$\!\pm\!$6.60e-4 & \bft{2.37e-2}$\!\pm\!$\bft{1.10e-2} & \bft{1.15e-2}$\!\pm\!$\bft{4.70e-3} & \bft{1.18e-2}$\!\pm\!$\bft{5.54e-3} & \bft{6.60e-3}$\!\pm\!$\bft{2.49e-3} \\
IL5 & 1.35e-3$\!\pm\!$3.25e-4 & 1.62e-3$\!\pm\!$6.56e-4 & 2.38e-2$\!\pm\!$1.10e-2 & 1.16e-2$\!\pm\!$4.67e-3 & 1.19e-2$\!\pm\!$5.53e-3 & 6.64e-3$\!\pm\!$2.47e-3 \\
IL9 & \bft{1.33e-3}$\!\pm\!$\bft{3.18e-4} & \bft{1.62e-3}$\!\pm\!$\bft{6.51e-4} & 2.39e-2$\!\pm\!$1.09e-2 & 1.17e-2$\!\pm\!$4.61e-3 & 1.20e-2$\!\pm\!$5.48e-3 & 6.67e-3$\!\pm\!$2.43e-3 \\
PS  & 1.57e-3$\!\pm\!$3.24e-4 & 1.66e-3$\!\pm\!$6.38e-4 & 2.39e-2$\!\pm\!$1.08e-2 & 1.18e-2$\!\pm\!$4.52e-3 & 1.20e-2$\!\pm\!$5.45e-3 & 6.83e-3$\!\pm\!$2.34e-3 \\
SF  & 3.86e-3$\!\pm\!$2.01e-3 & 3.69e-3$\!\pm\!$3.15e-3 & 4.33e-2$\!\pm\!$1.82e-2 & 3.73e-2$\!\pm\!$9.16e-3 & 2.17e-2$\!\pm\!$9.21e-3 & 1.95e-2$\!\pm\!$4.32e-3 \\
\bottomrule
\end{tabular}}
\end{table*}

\paragraph{Key findings}
Across the six physical magnitudes reported in Table~\ref{tab:summary_case30}, PS reduces the RMSE by 45--68\% relative to SF, with the largest gains concentrated in the reactive-power injections and flows. This outcome matches the structural expectation that enforcing the AC power flow equations through the WLS solve at inference yields estimates more faithful to the physical grid than the unconstrained mapping learned by~SF.
%
Regarding the IL configurations, a clear trade-off is observed: increasing $\gamma$ refines the accuracy of state variables but leads to higher errors in power derived magnitudes. Consequently, IL9 yields the highest precision for voltage magnitudes and angles, whereas IL1 excels in estimating power flows and injections. Crucially, the proposed IL framework consistently outperforms both the SF and PS baselines across all six physical magnitudes. Quantitatively, the best-performing IL configuration reduces RMSE by 0.8--15\% relative to PS and by 45--69\% relative to SF across the six magnitudes. The improvement is largest for voltage magnitude (15\% over PS) and smallest for active-power injection (under 1\% over PS), confirming that the gains from end-to-end training are concentrated in the state variables that the WLS sensitivity~\eqref{eq:wls_jacobian} propagates most directly.

\paragraph{Sensitivity to experimental design parameters}
Analysis of the disaggregated results in Table \ref{tab:results_case_30} (\ref{appendix:detailed_results}) confirms that the observed trends are consistent across different network conditions, highlighting the framework's generality. Two primary factors dictate the performance across all tested methodologies. First, demand variability acts as a primary driver of estimation error: as variability increases, the RMSE rises across all configurations because the NN must generalize over a more complex and diverse operating manifold. Second, the density of real-time measurements obviously influences accuracy, where higher sensing density of PMUs provides more direct observations on the whole state of the network. While these external factors impact all approaches, the relative performance ranking remains stable, with the IL framework maintaining a lead over the SF and PS benchmarks across the entire range of variability and measurement scenarios.

\subsubsection{IEEE 33-bus passive DN}


\begin{table*}[!t]
\caption{Average results for IEEE 33-bus (mean $\pm$ std) for different variability levels, measurement configurations, and data training sets.}
\label{tab:summary_case33}
\centering
\scriptsize
\renewcommand{\arraystretch}{1.0}
\setlength{\tabcolsep}{2pt}
\newcommand{\na}{---}
\resizebox{\textwidth}{!}{
\renewcommand{\arraystretch}{1.2}
\begin{tabular}{ccccccc}
\toprule
 & \textbf{RMSE$_{V}$} & \textbf{RMSE$_{\theta}$} & \textbf{RMSE$_{P}$} & \textbf{RMSE$_{Q}$} & \textbf{RMSE$_{P_f}$} & \textbf{RMSE$_{Q_f}$} \\
\midrule
IL1 & 1.67e-3$\!\pm\!$6.65e-4 & 2.16e-3$\!\pm\!$9.85e-4 & \bft{2.38e-2}$\!\pm\!$\bft{1.83e-2} & \bft{2.86e-2}$\!\pm\!$\bft{2.16e-2} & \bft{1.76e-2}$\!\pm\!$\bft{1.23e-2} & \bft{1.94e-2}$\!\pm\!$\bft{1.35e-2} \\
IL5 & \bft{1.66e-3}$\!\pm\!$\bft{6.31e-4} & 2.12e-3$\!\pm\!$9.32e-4 & 2.45e-2$\!\pm\!$1.87e-2 & 2.93e-2$\!\pm\!$2.20e-2 & 1.80e-2$\!\pm\!$1.26e-2 & 1.98e-2$\!\pm\!$1.38e-2 \\
IL9 & 1.68e-3$\!\pm\!$6.35e-4 & 2.11e-3$\!\pm\!$9.13e-4 & 2.49e-2$\!\pm\!$1.92e-2 & 2.96e-2$\!\pm\!$2.24e-2 & 1.84e-2$\!\pm\!$1.29e-2 & 2.01e-2$\!\pm\!$1.40e-2 \\
PS  & 1.94e-3$\!\pm\!$8.14e-4 & 2.68e-3$\!\pm\!$1.12e-3 & 2.42e-2$\!\pm\!$1.83e-2 & 2.90e-2$\!\pm\!$2.16e-2 & 1.81e-2$\!\pm\!$1.22e-2 & 1.99e-2$\!\pm\!$1.34e-2 \\
SF  & 2.46e-3$\!\pm\!$6.26e-4 & \bft{2.07e-3}$\!\pm\!$\bft{1.01e-3} & 7.42e-2$\!\pm\!$1.56e-2 & 7.14e-2$\!\pm\!$1.82e-2 & 4.84e-2$\!\pm\!$1.01e-2 & 4.36e-2$\!\pm\!$1.05e-2 \\
\bottomrule
\end{tabular}}
\end{table*}

\paragraph{Key findings}
Similarly, Table \ref{tab:summary_case33} presents the average RMSE results for the 33-bus radial DN. The relative performance of the best IL reduces RMSE by 1.4--21\% relative to PS and by 33--68\% relative to SF in five of the six magnitudes. 
The performance across methodologies reinforces the trade-offs observed in the meshed case: IL9 generally minimizes state variable errors, whereas IL1 provides superior accuracy for derived magnitudes. Notably, in this simpler radial topology, even though SF achieves in average the lowest RMSE for voltage angles, these estimation gains are misleading as they do not translate to reliable physical quantities. In fact, its corresponding RMSE for power injections and flows is, on average, 2 to 3 times higher than that of the proposed IL. 

\paragraph{Sensitivity to experimental design parameters}
The disaggregated results for the 33-bus radial network, detailed in Table \ref{tab:results_case33} (\ref{appendix:detailed_results}), reveal distinct topological sensitivities. Similarly to the 30-bus ADN, increasing levels of demand variability predictably increase the RMSE across all metrics. 
The effect of measurement configuration is, however, qualitatively different here, 
as BIF, END, and PMU differ not only in the density of real-time sensors but also in the type and placement of the available information: power-flow readings at bifurcation branches, voltage magnitudes at leaf nodes, and a single $V,\theta$ PMU at a strategic midpoint, respectively. 
While both the absolute error level and the magnitude of the IL-over-PS improvement vary across these scenarios, the relative ranking is preserved throughout the variability/configuration grid, with IL retaining its lead at 20\% variability and confirming that the framework's benefit does not rely on operating conditions being close to nominal.

\subsubsection{Voltage magnitude performance}

To conclude these case studies, Table \ref{tab:unified_improvement} summarizes the percentage improvement in voltage magnitude RMSE achieved by the best-performing IL approach relative to the PS baseline for both the meshed 30-bus and radial 33-bus DNs. We select this representative magnitude as it is crucial to ensure a stable operation and prevent equipment damage.
Across the 18 variability-configuration combinations, IL improves RMSE$_V$ over PS in every cell, with gains ranging from 2.9\% (30-bus, LOW, 20\% variability) to 33.7\% (30-bus, HIG, 5\% variability). Two patterns stand out. First, the IL gain is  smaller at highest variability than at lower variability levels, especially on the 30-bus network: averaged across configurations, it is 22\% at 5\%, 20\% at 10\%, and 6\% at 20\%; on the 33-bus the corresponding averages are 17\%, 17\%, and 10\%, showing a similar but more contained drop. End-to-end training therefore delivers its largest benefits when the system operates closer to nominal conditions, although it remains beneficial under high variability. Second, the configuration where the IL gain is less pronounced is LOW for the 30-bus system (the fewest PMUs) and END for the 33-bus one, the only configuration that measures voltage magnitudes directly at multiple buses and therefore leaves the smallest room for improvement on RMSE$_V$ specifically. Taken together, the table confirms that IL is never inferior to PS on voltage estimation and provides substantial improvements in most regimes, with the magnitude of the benefit contingent on both the variability of the operating point and the structure of the available real-time measurements. 

\begin{table}[h!]
\centering
\caption{Relative RMSE$_{V}$ improvement (\%) of the best IL approach over the PS baseline across both test systems, for different demand variability levels (Var.) and the different measurement configurations defined in Table~\ref{tab:scenarios}.}
\label{tab:unified_improvement}
\renewcommand{\arraystretch}{1.2}
\begin{tabular}{@{}c ccc p{0.1cm} ccc@{}}
\toprule
\multirow{2.5}{*}{\textbf{Var. (\%)}} & \multicolumn{3}{c}{\textbf{IEEE 30-bus System}} & & \multicolumn{3}{c}{\textbf{IEEE 33-bus System}} \\
\cmidrule(lr){2-4} \cmidrule(lr){6-8}
 & {LOW} & {MED} & {HIG} & & {BIF} & {END} & {PMU} \\ \midrule
5  & 15.5\% & 17.6\% & 33.7\% & & 24.4\% & 5.2\%  & 21.1\% \\
10 & 5.2\%  & 28.1\% & 25.5\% & & 9.1\%  & 12.1\% & 30.7\% \\
20 & 2.9\%  & 9.3\%  & 6.2\%  & & 13.7\% & 3.4\%  & 11.5\% \\ \bottomrule
\end{tabular}
\end{table}

\section{Conclusions}\label{sec:concl}

This work proposes an end-to-end learning pseudo-measurement generation methodology for both active and passive DNs with limited observability, based on 
embedding a WLS state estimator as an implicit differentiable layer inside a NN.
By integrating physical constraints during training, the proposed approach bridges the gap between direct forecasting and classical SE.

The results obtained for both the IEEE 30-bus and IEEE 33-bus systems demonstrate that the implicit-layer approaches consistently outperform the baseline methods under different variability levels and sensor placement configurations in terms of RMSE. 
While the state forecasting strategy may achieve low errors in some state variables, it often produces physically inconsistent solutions with large errors in power injections and flows. 
In contrast, the proposed implicit-layer framework provides more accurate and physically coherent estimates, achieving improvements of up to $33.7\%$ in voltage magnitude RMSE with respect to the pseudo-measurement generation baseline. These results demonstrate that integrating physical grid constraints directly into the pseudo-measurement reconstruction process through an implicit layer acts as a powerful corrective mechanism, enabling the state estimator to attain a level of physical consistency and estimation accuracy beyond what traditional decoupled or purely data-driven forecasting approaches can achieve.



\section*{Acknowledgment}
This work was supported by the Spanish Ministry of Science, Innovation and Universities (AEI/10.13039/501100011033) through project PID2023-148291NB-I00 and by the Department of Universities, Research and Innovation of the Regional Government of Andalusia through project DGP\_PIDI\_2024\_00851. The work of J.~G.~De la Varga was supported by the Spanish Ministry of Science, Innovation and Universities training program for PhDs with fellowship number PRE2021-098958. The authors thankfully acknowledge the computer resources (Picasso Supercomputer), technical expertise and assistance provided by the SCBI (Supercomputing and Bioinformatics) center of the University of M\'alaga.


\appendix
\section{Measurement function formulation}
\label{appendix:h}

In this work, a fixed topology is considered, where the measurement functions $\mathbf{h}$ relate the set of measurable outputs $\mathbf{z}$ to the system's state variables $(V_i,\theta_i)$. With a slight abuse of notation, this relationship is formulated as follows:
\begin{subequations} \label{eq:h}
\begin{align} 
    &h_{V_i} = V_i \\
    &h_{\theta_i} = \theta_i \\    
    &h_{P_{ij}} = V_iV_j \left( G_{ij} \cos{\theta_{ij}} + B_{ij} \sin{\theta_{ij}} \right) - G_{ij}V_i^2 \\    
    &h_{Q_{ij}} = V_iV_j \left( G_{ij} \sin{\theta_{ij}} - B_{ij} \cos{\theta_{ij}} \right) + V_i^2\left(B_{ij}-b_{ij}^S/2\right)\\
    &h_{P_i} = V_i \sum _{j} V_j \left( G_{ij} \cos{\theta_{ij}} + B_{ij} \sin{\theta_{ij}} \right) \\
    &h_{Q_i} = V_i \sum _{j} V_j \left( G_{ij} \sin{\theta_{ij}} - B_{ij} \cos{\theta_{ij}} \right)
\end{align}
\end{subequations}
where $G_{ij}$ and $B_{ij}$ are, respectively, the real and imaginary parts of the bus admittance matrix, $b_{ij}^S$ is the shunt susceptance of the line between nodes $i\text{ and }j$, and ${\theta_{ij}=\theta_i - \theta_j}$.

\section{Detailed results}
\label{appendix:detailed_results}

Tables~\ref{tab:results_case_30} and~\ref{tab:results_case33} report the full per-configuration results of the three data-driven approaches SF, PS and IL; for the IEEE 30-bus and IEEE 33-bus networks, respectively. Each table details the average and standard deviation over seeds across the three demand variability levels and the different measurement configurations (detailed in Table~\ref{tab:scenarios}), complementing the aggregated summary in Section~\ref{sec:results}.

\begin{table}[!ht]
\caption{Detailed results for IEEE 30-bus (mean $\pm$ std over seeds) for different demand variability levels (Var.), measurement configurations (Conf.) and all the different approaches (App.).}
\label{tab:results_case_30}
\centering
\resizebox{\textwidth}{!}{
\begin{tabular}{ccccccccc}
\toprule
\bft{Var.} & \bft{Conf.} & \bft{App.} & \bft{RMSE$_{V}$} & \bft{RMSE$_{\theta}$} & \bft{RMSE$_{P}$} & \bft{RMSE$_{Q}$} & \bft{RMSE$_{P_f}$} & \bft{RMSE$_{Q_f}$} \\
\midrule
\multirow{15}{*}{5\%} & \multirow{5}{*}{HIG} & IL1 & 8.87e-4$\!\pm\!$2.67e-5 & \bft{8.27e-4}$\!\pm\!$\bft{1.84e-5} & \bft{1.14e-2}$\!\pm\!$\bft{3.28e-4} & \bft{6.00e-3}$\!\pm\!$\bft{1.80e-4} & \bft{5.69e-3}$\!\pm\!$\bft{1.13e-4} & \bft{3.47e-3}$\!\pm\!$\bft{1.07e-4} \\
 &  & IL5 & \bft{8.62e-4}$\!\pm\!$\bft{5.23e-5} & 8.31e-4$\!\pm\!$1.75e-5 & 1.17e-2$\!\pm\!$3.33e-4 & 6.15e-3$\!\pm\!$1.77e-4 & 5.80e-3$\!\pm\!$1.21e-4 & 3.53e-3$\!\pm\!$1.15e-4 \\
 &  & IL9 & 8.74e-4$\!\pm\!$1.73e-5 & 8.41e-4$\!\pm\!$1.29e-5 & 1.19e-2$\!\pm\!$3.05e-4 & 6.32e-3$\!\pm\!$1.90e-4 & 5.93e-3$\!\pm\!$9.54e-5 & 3.61e-3$\!\pm\!$9.87e-5 \\
 &  & PS & 1.30e-3$\!\pm\!$4.94e-5 & 8.97e-4$\!\pm\!$2.62e-5 & 1.19e-2$\!\pm\!$2.88e-4 & 6.69e-3$\!\pm\!$1.20e-4 & 5.94e-3$\!\pm\!$7.81e-5 & 3.98e-3$\!\pm\!$1.04e-4 \\
 &  & SF & 6.38e-3$\!\pm\!$7.69e-5 & 3.60e-3$\!\pm\!$7.85e-5 & 2.62e-2$\!\pm\!$4.43e-4 & 3.02e-2$\!\pm\!$1.86e-3 & 1.39e-2$\!\pm\!$1.46e-4 & 1.93e-2$\!\pm\!$5.98e-4 \\
\cmidrule{2-9}
 & \multirow{5}{*}{LOW} & IL1 & 1.11e-3$\!\pm\!$4.00e-5 & 1.20e-3$\!\pm\!$2.27e-5 & \bft{1.45e-2}$\!\pm\!$\bft{1.68e-4} & \bft{6.97e-3}$\!\pm\!$\bft{1.09e-4} & 7.25e-3$\!\pm\!$4.32e-5 & \bft{4.20e-3}$\!\pm\!$\bft{7.84e-5} \\
 &  & IL5 & 1.08e-3$\!\pm\!$5.25e-5 & \bft{1.19e-3}$\!\pm\!$\bft{1.32e-5} & 1.45e-2$\!\pm\!$2.91e-4 & 6.99e-3$\!\pm\!$2.46e-4 & \bft{7.23e-3}$\!\pm\!$\bft{1.28e-4} & 4.21e-3$\!\pm\!$1.56e-4 \\
 &  & IL9 & \bft{1.08e-3}$\!\pm\!$\bft{4.47e-5} & 1.20e-3$\!\pm\!$1.35e-5 & 1.47e-2$\!\pm\!$2.22e-4 & 7.18e-3$\!\pm\!$1.51e-4 & 7.32e-3$\!\pm\!$7.83e-5 & 4.29e-3$\!\pm\!$9.89e-5 \\
 &  & PS & 1.28e-3$\!\pm\!$1.45e-4 & 1.24e-3$\!\pm\!$2.04e-5 & 1.49e-2$\!\pm\!$2.14e-4 & 7.25e-3$\!\pm\!$2.57e-4 & 7.42e-3$\!\pm\!$8.92e-5 & 4.46e-3$\!\pm\!$2.01e-4 \\
 &  & SF & 8.32e-3$\!\pm\!$1.18e-4 & 1.20e-2$\!\pm\!$4.33e-4 & 8.29e-2$\!\pm\!$3.61e-3 & 4.34e-2$\!\pm\!$3.73e-3 & 4.16e-2$\!\pm\!$1.65e-3 & 2.58e-2$\!\pm\!$1.38e-3 \\
\cmidrule{2-9}
 & \multirow{5}{*}{MED} & IL1 & 1.54e-3$\!\pm\!$1.76e-4 & \bft{1.15e-3}$\!\pm\!$\bft{6.28e-5} & \bft{1.27e-2}$\!\pm\!$\bft{5.22e-4} & \bft{6.94e-3}$\!\pm\!$\bft{5.05e-4} & \bft{6.25e-3}$\!\pm\!$\bft{2.12e-4} & \bft{4.53e-3}$\!\pm\!$\bft{3.79e-4} \\
 &  & IL5 & 1.54e-3$\!\pm\!$1.17e-4 & 1.18e-3$\!\pm\!$4.34e-5 & 1.28e-2$\!\pm\!$3.85e-4 & 7.08e-3$\!\pm\!$3.78e-4 & 6.33e-3$\!\pm\!$1.39e-4 & 4.58e-3$\!\pm\!$2.48e-4 \\
 &  & IL9 & \bft{1.53e-3}$\!\pm\!$\bft{1.57e-4} & 1.18e-3$\!\pm\!$4.21e-5 & 1.31e-2$\!\pm\!$3.46e-4 & 7.28e-3$\!\pm\!$3.98e-4 & 6.44e-3$\!\pm\!$1.27e-4 & 4.68e-3$\!\pm\!$3.08e-4 \\
 &  & PS & 1.86e-3$\!\pm\!$7.19e-5 & 1.26e-3$\!\pm\!$1.39e-5 & 1.30e-2$\!\pm\!$3.87e-4 & 7.34e-3$\!\pm\!$3.84e-4 & 6.42e-3$\!\pm\!$1.32e-4 & 4.90e-3$\!\pm\!$2.54e-4 \\
 &  & SF & 2.92e-3$\!\pm\!$1.18e-4 & 1.58e-3$\!\pm\!$3.47e-5 & 2.37e-2$\!\pm\!$1.10e-3 & 2.45e-2$\!\pm\!$2.11e-3 & 1.14e-2$\!\pm\!$5.47e-4 & 1.25e-2$\!\pm\!$9.52e-4 \\
\midrule
\multirow{15}{*}{10\%} & \multirow{5}{*}{HIG} & IL1 & 1.14e-3$\!\pm\!$4.44e-5 & 1.10e-3$\!\pm\!$1.37e-5 & \bft{1.85e-2}$\!\pm\!$\bft{2.48e-4} & \bft{9.67e-3}$\!\pm\!$\bft{1.51e-4} & \bft{9.29e-3}$\!\pm\!$\bft{1.07e-4} & \bft{5.42e-3}$\!\pm\!$\bft{8.41e-5} \\
 &  & IL5 & \bft{1.04e-3}$\!\pm\!$\bft{2.63e-5} & \bft{1.10e-3}$\!\pm\!$\bft{1.00e-5} & 1.86e-2$\!\pm\!$2.20e-4 & 9.77e-3$\!\pm\!$1.74e-4 & 9.34e-3$\!\pm\!$9.70e-5 & 5.45e-3$\!\pm\!$8.47e-5 \\
 &  & IL9 & 1.04e-3$\!\pm\!$2.75e-5 & 1.10e-3$\!\pm\!$1.07e-5 & 1.88e-2$\!\pm\!$3.00e-4 & 9.85e-3$\!\pm\!$2.04e-4 & 9.40e-3$\!\pm\!$1.36e-4 & 5.48e-3$\!\pm\!$1.01e-4 \\
 &  & PS & 1.39e-3$\!\pm\!$8.50e-5 & 1.14e-3$\!\pm\!$1.82e-5 & 1.88e-2$\!\pm\!$2.46e-4 & 1.01e-2$\!\pm\!$1.97e-4 & 9.44e-3$\!\pm\!$1.06e-4 & 5.71e-3$\!\pm\!$1.19e-4 \\
 &  & SF & 2.66e-3$\!\pm\!$7.46e-5 & 1.81e-3$\!\pm\!$8.24e-5 & 2.90e-2$\!\pm\!$2.20e-3 & 3.98e-2$\!\pm\!$6.79e-3 & 1.39e-2$\!\pm\!$7.37e-4 & 1.77e-2$\!\pm\!$2.69e-3 \\
\cmidrule{2-9}
 & \multirow{5}{*}{LOW} & IL1 & 1.40e-3$\!\pm\!$3.90e-5 & 1.85e-3$\!\pm\!$6.11e-5 & \bft{2.66e-2}$\!\pm\!$\bft{3.08e-4} & \bft{1.19e-2}$\!\pm\!$\bft{2.27e-4} & \bft{1.34e-2}$\!\pm\!$\bft{1.70e-4} & 6.91e-3$\!\pm\!$1.31e-4 \\
 &  & IL5 & 1.39e-3$\!\pm\!$1.75e-5 & 1.83e-3$\!\pm\!$2.66e-5 & 2.67e-2$\!\pm\!$2.56e-4 & 1.20e-2$\!\pm\!$1.03e-4 & 1.35e-2$\!\pm\!$1.26e-4 & 6.93e-3$\!\pm\!$4.65e-5 \\
 &  & IL9 & \bft{1.37e-3}$\!\pm\!$\bft{6.33e-5} & \bft{1.82e-3}$\!\pm\!$\bft{4.11e-5} & 2.67e-2$\!\pm\!$3.34e-4 & 1.20e-2$\!\pm\!$1.43e-4 & 1.35e-2$\!\pm\!$1.65e-4 & \bft{6.91e-3}$\!\pm\!$\bft{1.05e-4} \\
 &  & PS & 1.44e-3$\!\pm\!$1.87e-5 & 1.86e-3$\!\pm\!$2.45e-5 & 2.68e-2$\!\pm\!$1.98e-4 & 1.21e-2$\!\pm\!$1.14e-4 & 1.35e-2$\!\pm\!$7.69e-5 & 6.98e-3$\!\pm\!$5.60e-5 \\
 &  & SF & 2.68e-3$\!\pm\!$8.88e-5 & 2.11e-3$\!\pm\!$6.36e-5 & 4.16e-2$\!\pm\!$2.85e-3 & 4.77e-2$\!\pm\!$3.92e-3 & 2.11e-2$\!\pm\!$1.33e-3 & 2.35e-2$\!\pm\!$1.89e-3 \\
\cmidrule{2-9}
 & \multirow{5}{*}{MED} & IL1 & \bft{1.50e-3}$\!\pm\!$\bft{2.76e-5} & \bft{1.66e-3}$\!\pm\!$\bft{6.41e-6} & \bft{2.05e-2}$\!\pm\!$\bft{2.57e-4} & \bft{1.05e-2}$\!\pm\!$\bft{6.81e-5} & \bft{1.02e-2}$\!\pm\!$\bft{7.89e-5} & \bft{6.17e-3}$\!\pm\!$\bft{3.11e-5} \\
 &  & IL5 & 1.65e-3$\!\pm\!$9.00e-5 & 1.67e-3$\!\pm\!$9.20e-6 & 2.08e-2$\!\pm\!$2.06e-4 & 1.08e-2$\!\pm\!$1.06e-4 & 1.03e-2$\!\pm\!$6.66e-5 & 6.34e-3$\!\pm\!$8.86e-5 \\
 &  & IL9 & 1.55e-3$\!\pm\!$4.29e-5 & 1.67e-3$\!\pm\!$5.85e-6 & 2.08e-2$\!\pm\!$1.99e-4 & 1.08e-2$\!\pm\!$3.55e-5 & 1.03e-2$\!\pm\!$5.07e-5 & 6.30e-3$\!\pm\!$8.95e-6 \\
 &  & PS & 2.09e-3$\!\pm\!$1.85e-4 & 1.72e-3$\!\pm\!$9.72e-6 & 2.09e-2$\!\pm\!$2.11e-4 & 1.10e-2$\!\pm\!$1.27e-4 & 1.04e-2$\!\pm\!$6.90e-5 & 6.63e-3$\!\pm\!$1.06e-4 \\
 &  & SF & 2.96e-3$\!\pm\!$1.52e-4 & 2.11e-3$\!\pm\!$1.42e-5 & 3.06e-2$\!\pm\!$1.17e-3 & 2.79e-2$\!\pm\!$1.30e-3 & 1.49e-2$\!\pm\!$4.80e-4 & 1.46e-2$\!\pm\!$1.42e-3 \\
\midrule
\multirow{15}{*}{20\%} & \multirow{5}{*}{HIG} & IL1 & 1.33e-3$\!\pm\!$5.23e-5 & 1.44e-3$\!\pm\!$2.88e-5 & \bft{2.69e-2}$\!\pm\!$\bft{2.05e-4} & \bft{1.41e-2}$\!\pm\!$\bft{1.07e-4} & \bft{1.35e-2}$\!\pm\!$\bft{1.27e-4} & \bft{7.83e-3}$\!\pm\!$\bft{5.35e-5} \\
 &  & IL5 & 1.31e-3$\!\pm\!$4.88e-5 & 1.44e-3$\!\pm\!$2.77e-5 & 2.70e-2$\!\pm\!$2.34e-4 & 1.42e-2$\!\pm\!$1.45e-4 & 1.35e-2$\!\pm\!$1.55e-4 & 7.85e-3$\!\pm\!$8.38e-5 \\
 &  & IL9 & \bft{1.27e-3}$\!\pm\!$\bft{4.61e-5} & \bft{1.43e-3}$\!\pm\!$\bft{2.90e-5} & 2.72e-2$\!\pm\!$3.56e-4 & 1.42e-2$\!\pm\!$2.23e-4 & 1.36e-2$\!\pm\!$2.15e-4 & 7.87e-3$\!\pm\!$1.27e-4 \\
 &  & PS & 1.35e-3$\!\pm\!$5.06e-5 & 1.45e-3$\!\pm\!$3.88e-5 & 2.71e-2$\!\pm\!$3.34e-4 & 1.43e-2$\!\pm\!$2.19e-4 & 1.36e-2$\!\pm\!$2.14e-4 & 7.91e-3$\!\pm\!$1.19e-4 \\
 &  & SF & 2.81e-3$\!\pm\!$5.46e-5 & 2.66e-3$\!\pm\!$4.03e-5 & 4.76e-2$\!\pm\!$2.60e-3 & 4.69e-2$\!\pm\!$9.55e-3 & 2.37e-2$\!\pm\!$1.09e-3 & 2.17e-2$\!\pm\!$3.34e-3 \\
\cmidrule{2-9}
 & \multirow{5}{*}{LOW} & IL1 & 1.50e-3$\!\pm\!$4.67e-5 & 2.73e-3$\!\pm\!$4.17e-5 & 4.56e-2$\!\pm\!$5.87e-4 & \bft{1.95e-2}$\!\pm\!$\bft{2.79e-4} & 2.29e-2$\!\pm\!$3.16e-4 & 1.08e-2$\!\pm\!$1.57e-4 \\
 &  & IL5 & 1.45e-3$\!\pm\!$3.92e-5 & 2.73e-3$\!\pm\!$5.18e-5 & 4.57e-2$\!\pm\!$7.86e-4 & 1.96e-2$\!\pm\!$3.10e-4 & 2.29e-2$\!\pm\!$4.28e-4 & 1.08e-2$\!\pm\!$1.74e-4 \\
 &  & IL9 & \bft{1.45e-3}$\!\pm\!$\bft{3.12e-5} & \bft{2.72e-3}$\!\pm\!$\bft{3.95e-5} & 4.56e-2$\!\pm\!$5.99e-4 & 1.96e-2$\!\pm\!$2.54e-4 & 2.29e-2$\!\pm\!$3.26e-4 & \bft{1.08e-2}$\!\pm\!$\bft{1.49e-4} \\
 &  & PS & 1.49e-3$\!\pm\!$5.24e-5 & 2.73e-3$\!\pm\!$3.92e-5 & \bft{4.55e-2}$\!\pm\!$\bft{5.92e-4} & 1.95e-2$\!\pm\!$2.83e-4 & \bft{2.28e-2}$\!\pm\!$\bft{3.17e-4} & 1.08e-2$\!\pm\!$1.60e-4 \\
 &  & SF & 2.70e-3$\!\pm\!$7.82e-5 & 3.15e-3$\!\pm\!$4.87e-5 & 5.67e-2$\!\pm\!$1.26e-3 & 4.07e-2$\!\pm\!$2.00e-3 & 2.88e-2$\!\pm\!$6.20e-4 & 2.12e-2$\!\pm\!$9.80e-4 \\
\cmidrule{2-9}
 & \multirow{5}{*}{MED} & IL1 & 2.03e-3$\!\pm\!$8.29e-5 & 2.64e-3$\!\pm\!$2.52e-5 & \bft{3.50e-2}$\!\pm\!$\bft{9.84e-5} & 1.77e-2$\!\pm\!$1.71e-4 & 1.75e-2$\!\pm\!$9.38e-5 & 9.95e-3$\!\pm\!$1.39e-4 \\
 &  & IL5 & 1.94e-3$\!\pm\!$5.83e-5 & 2.63e-3$\!\pm\!$2.89e-5 & 3.51e-2$\!\pm\!$1.00e-4 & 1.77e-2$\!\pm\!$1.46e-4 & \bft{1.75e-2}$\!\pm\!$\bft{5.62e-5} & \bft{9.92e-3}$\!\pm\!$\bft{1.13e-4} \\
 &  & IL9 & \bft{1.94e-3}$\!\pm\!$\bft{3.46e-5} & \bft{2.63e-3}$\!\pm\!$\bft{2.50e-5} & 3.52e-2$\!\pm\!$9.79e-5 & \bft{1.77e-2}$\!\pm\!$\bft{1.33e-4} & 1.76e-2$\!\pm\!$6.69e-5 & 9.94e-3$\!\pm\!$9.79e-5 \\
 &  & PS & 2.14e-3$\!\pm\!$9.47e-5 & 2.66e-3$\!\pm\!$2.49e-5 & 3.51e-2$\!\pm\!$4.78e-5 & 1.78e-2$\!\pm\!$1.30e-4 & 1.76e-2$\!\pm\!$4.33e-5 & 1.00e-2$\!\pm\!$1.10e-4 \\
 &  & SF & 2.93e-3$\!\pm\!$9.05e-5 & 3.47e-3$\!\pm\!$6.15e-5 & 4.64e-2$\!\pm\!$5.89e-4 & 3.11e-2$\!\pm\!$1.49e-3 & 2.30e-2$\!\pm\!$2.41e-4 & 1.72e-2$\!\pm\!$7.99e-4 \\
\bottomrule
\end{tabular}
}
\end{table}

\begin{table*}[!ht]
\caption{Detailed results for IEEE 33-bus (mean $\pm$ std over seeds) for different demand variability levels (Var.), measurement configurations (Conf.) and all the different approaches (App.).}
\label{tab:results_case33}
\centering
\resizebox{\textwidth}{!}{
\begin{tabular}{ccccccccc}
\toprule
\bft{Var.} & \bft{Conf.} & \bft{App.} & \bft{RMSE$_{V}$} & \bft{RMSE$_{\theta}$} & \bft{RMSE$_{P}$} & \bft{RMSE$_{Q}$} & \bft{RMSE$_{P_f}$} & \bft{RMSE$_{Q_f}$} \\
\midrule
\multirow{15}{*}{5\%} & \multirow{5}{*}{BIF} & IL1 & \bft{1.95e-3}$\!\pm\!$\bft{3.55e-5} & 8.79e-4$\!\pm\!$5.89e-5 & \bft{3.60e-3}$\!\pm\!$\bft{9.71e-5} & \bft{4.86e-3}$\!\pm\!$\bft{1.37e-4} & \bft{3.95e-3}$\!\pm\!$\bft{1.36e-4} & \bft{4.28e-3}$\!\pm\!$\bft{1.01e-4} \\
 &  & IL5 & 2.04e-3$\!\pm\!$1.13e-4 & \bft{8.40e-4}$\!\pm\!$\bft{5.00e-5} & 3.84e-3$\!\pm\!$1.72e-4 & 5.11e-3$\!\pm\!$2.13e-4 & 4.17e-3$\!\pm\!$2.24e-4 & 4.48e-3$\!\pm\!$2.07e-4 \\
 &  & IL9 & 2.06e-3$\!\pm\!$6.99e-5 & 8.93e-4$\!\pm\!$3.72e-5 & 4.11e-3$\!\pm\!$1.69e-4 & 5.34e-3$\!\pm\!$1.86e-4 & 4.38e-3$\!\pm\!$2.10e-4 & 4.69e-3$\!\pm\!$9.65e-5 \\
 &  & PS & 2.58e-3$\!\pm\!$3.44e-4 & 1.27e-3$\!\pm\!$1.33e-4 & 4.07e-3$\!\pm\!$1.08e-4 & 5.33e-3$\!\pm\!$1.95e-4 & 4.93e-3$\!\pm\!$4.58e-4 & 5.00e-3$\!\pm\!$1.60e-4 \\
 &  & SF & 2.50e-3$\!\pm\!$7.62e-5 & 1.14e-3$\!\pm\!$6.10e-5 & 7.18e-2$\!\pm\!$6.59e-3 & 7.70e-2$\!\pm\!$4.35e-3 & 4.33e-2$\!\pm\!$3.15e-3 & 4.43e-2$\!\pm\!$2.39e-3 \\
\cmidrule{2-9}
 & \multirow{5}{*}{END} & IL1 & \bft{8.22e-4}$\!\pm\!$\bft{1.85e-5} & 1.42e-3$\!\pm\!$3.80e-4 & \bft{1.62e-2}$\!\pm\!$\bft{2.05e-4} & \bft{1.89e-2}$\!\pm\!$\bft{1.69e-4} & \bft{1.17e-2}$\!\pm\!$\bft{9.42e-5} & \bft{1.27e-2}$\!\pm\!$\bft{5.78e-5} \\
 &  & IL5 & 8.41e-4$\!\pm\!$3.41e-5 & 1.48e-3$\!\pm\!$2.21e-4 & 1.63e-2$\!\pm\!$2.33e-4 & 1.91e-2$\!\pm\!$1.84e-4 & 1.19e-2$\!\pm\!$2.11e-4 & 1.28e-2$\!\pm\!$1.77e-4 \\
 &  & IL9 & 9.49e-4$\!\pm\!$5.52e-5 & 1.45e-3$\!\pm\!$1.62e-4 & 1.68e-2$\!\pm\!$4.21e-4 & 1.95e-2$\!\pm\!$4.39e-4 & 1.23e-2$\!\pm\!$4.65e-4 & 1.33e-2$\!\pm\!$4.66e-4 \\
 &  & PS & 8.68e-4$\!\pm\!$6.40e-5 & 1.90e-3$\!\pm\!$1.98e-4 & 1.63e-2$\!\pm\!$2.04e-4 & 1.90e-2$\!\pm\!$1.72e-4 & 1.19e-2$\!\pm\!$1.27e-4 & 1.29e-2$\!\pm\!$1.01e-4 \\
 &  & SF & 1.60e-3$\!\pm\!$1.34e-4 & \bft{1.04e-3}$\!\pm\!$\bft{5.28e-6} & 8.70e-2$\!\pm\!$1.79e-2 & 5.76e-2$\!\pm\!$1.05e-2 & 5.50e-2$\!\pm\!$1.06e-2 & 3.47e-2$\!\pm\!$5.98e-3 \\
\cmidrule{2-9}
 & \multirow{5}{*}{PMU} & IL1 & 1.28e-3$\!\pm\!$2.61e-5 & 1.64e-3$\!\pm\!$2.89e-5 & \bft{1.61e-2}$\!\pm\!$\bft{1.84e-4} & \bft{1.89e-2}$\!\pm\!$\bft{2.25e-4} & \bft{1.21e-2}$\!\pm\!$\bft{1.19e-4} & 1.28e-2$\!\pm\!$1.66e-4 \\
 &  & IL5 & 1.27e-3$\!\pm\!$4.85e-5 & 1.63e-3$\!\pm\!$2.15e-5 & 1.62e-2$\!\pm\!$1.97e-4 & 1.89e-2$\!\pm\!$2.33e-4 & 1.21e-2$\!\pm\!$1.36e-4 & \bft{1.28e-2}$\!\pm\!$\bft{1.50e-4} \\
 &  & IL9 & \bft{1.26e-3}$\!\pm\!$\bft{2.16e-5} & 1.63e-3$\!\pm\!$2.97e-5 & 1.63e-2$\!\pm\!$2.56e-4 & 1.90e-2$\!\pm\!$3.19e-4 & 1.22e-2$\!\pm\!$1.53e-4 & 1.30e-2$\!\pm\!$2.35e-4 \\
 &  & PS & 1.60e-3$\!\pm\!$1.58e-4 & 2.02e-3$\!\pm\!$1.73e-4 & 1.62e-2$\!\pm\!$1.93e-4 & 1.90e-2$\!\pm\!$2.33e-4 & 1.23e-2$\!\pm\!$1.03e-4 & 1.31e-2$\!\pm\!$2.24e-4 \\
 &  & SF & 1.73e-3$\!\pm\!$9.42e-5 & \bft{1.04e-3}$\!\pm\!$\bft{1.68e-5} & 4.87e-2$\!\pm\!$6.98e-3 & 4.10e-2$\!\pm\!$4.16e-3 & 3.14e-2$\!\pm\!$4.61e-3 & 2.55e-2$\!\pm\!$2.47e-3 \\
\midrule
\multirow{15}{*}{10\%} & \multirow{5}{*}{BIF} & IL1 & 2.33e-3$\!\pm\!$2.98e-5 & 1.41e-3$\!\pm\!$1.84e-4 & \bft{4.80e-3}$\!\pm\!$\bft{2.04e-5} & \bft{6.53e-3}$\!\pm\!$\bft{2.92e-5} & 5.13e-3$\!\pm\!$8.95e-5 & \bft{5.83e-3}$\!\pm\!$\bft{5.04e-5} \\
 &  & IL5 & \bft{2.23e-3}$\!\pm\!$\bft{5.69e-5} & 1.34e-3$\!\pm\!$5.58e-5 & 4.86e-3$\!\pm\!$7.77e-5 & 6.59e-3$\!\pm\!$6.90e-5 & \bft{5.09e-3}$\!\pm\!$\bft{9.08e-5} & 5.85e-3$\!\pm\!$7.80e-5 \\
 &  & IL9 & 2.28e-3$\!\pm\!$4.93e-5 & \bft{1.34e-3}$\!\pm\!$\bft{3.73e-5} & 5.41e-3$\!\pm\!$4.37e-4 & 6.96e-3$\!\pm\!$2.91e-4 & 5.56e-3$\!\pm\!$3.79e-4 & 6.21e-3$\!\pm\!$3.23e-4 \\
 &  & PS & 2.46e-3$\!\pm\!$1.19e-4 & 1.71e-3$\!\pm\!$3.65e-4 & 4.92e-3$\!\pm\!$5.99e-5 & 6.64e-3$\!\pm\!$3.05e-5 & 5.32e-3$\!\pm\!$1.71e-4 & 6.28e-3$\!\pm\!$3.05e-4 \\
 &  & SF & 3.02e-3$\!\pm\!$1.42e-4 & 1.58e-3$\!\pm\!$5.11e-5 & 6.91e-2$\!\pm\!$5.04e-3 & 1.03e-1$\!\pm\!$6.05e-3 & 4.55e-2$\!\pm\!$3.25e-3 & 5.87e-2$\!\pm\!$3.30e-3 \\
\cmidrule{2-9}
 & \multirow{5}{*}{END} & IL1 & \bft{8.97e-4}$\!\pm\!$\bft{1.77e-5} & 2.37e-3$\!\pm\!$3.11e-4 & \bft{2.99e-2}$\!\pm\!$\bft{2.11e-3} & \bft{3.55e-2}$\!\pm\!$\bft{2.16e-3} & \bft{2.12e-2}$\!\pm\!$\bft{1.47e-3} & \bft{2.34e-2}$\!\pm\!$\bft{1.42e-3} \\
 &  & IL5 & 1.01e-3$\!\pm\!$1.42e-4 & 2.63e-3$\!\pm\!$4.01e-4 & 3.30e-2$\!\pm\!$2.02e-3 & 3.86e-2$\!\pm\!$1.47e-3 & 2.35e-2$\!\pm\!$1.44e-3 & 2.56e-2$\!\pm\!$1.25e-3 \\
 &  & IL9 & 9.72e-4$\!\pm\!$9.91e-5 & 2.45e-3$\!\pm\!$2.23e-4 & 3.28e-2$\!\pm\!$1.29e-3 & 3.82e-2$\!\pm\!$8.41e-4 & 2.33e-2$\!\pm\!$9.97e-4 & 2.53e-2$\!\pm\!$6.14e-4 \\
 &  & PS & 1.02e-3$\!\pm\!$1.48e-4 & 3.57e-3$\!\pm\!$9.93e-4 & 3.21e-2$\!\pm\!$7.22e-4 & 3.78e-2$\!\pm\!$5.42e-4 & 2.29e-2$\!\pm\!$4.55e-4 & 2.53e-2$\!\pm\!$4.12e-4 \\
 &  & SF & 1.92e-3$\!\pm\!$1.53e-4 & \bft{1.94e-3}$\!\pm\!$\bft{2.05e-5} & 9.26e-2$\!\pm\!$1.77e-2 & 6.70e-2$\!\pm\!$9.47e-3 & 5.99e-2$\!\pm\!$1.12e-2 & 4.13e-2$\!\pm\!$5.41e-3 \\
\cmidrule{2-9}
 & \multirow{5}{*}{PMU} & IL1 & 1.63e-3$\!\pm\!$1.45e-4 & 2.17e-3$\!\pm\!$2.37e-4 & \bft{3.03e-2}$\!\pm\!$\bft{7.05e-4} & 3.58e-2$\!\pm\!$6.89e-4 & 2.18e-2$\!\pm\!$4.86e-4 & 2.35e-2$\!\pm\!$4.61e-4 \\
 &  & IL5 & 1.50e-3$\!\pm\!$2.60e-5 & 1.94e-3$\!\pm\!$5.27e-5 & 3.04e-2$\!\pm\!$6.82e-4 & 3.58e-2$\!\pm\!$6.89e-4 & 2.18e-2$\!\pm\!$4.94e-4 & 2.35e-2$\!\pm\!$4.71e-4 \\
 &  & IL9 & \bft{1.49e-3}$\!\pm\!$\bft{3.19e-5} & \bft{1.91e-3}$\!\pm\!$\bft{3.42e-5} & 3.04e-2$\!\pm\!$6.69e-4 & \bft{3.58e-2}$\!\pm\!$\bft{6.54e-4} & \bft{2.18e-2}$\!\pm\!$\bft{4.74e-4} & \bft{2.35e-2}$\!\pm\!$\bft{4.39e-4} \\
 &  & PS & 2.16e-3$\!\pm\!$5.87e-4 & 2.69e-3$\!\pm\!$4.56e-4 & 3.04e-2$\!\pm\!$6.85e-4 & 3.58e-2$\!\pm\!$6.68e-4 & 2.21e-2$\!\pm\!$5.72e-4 & 2.37e-2$\!\pm\!$4.53e-4 \\
 &  & SF & 2.21e-3$\!\pm\!$1.21e-4 & 1.96e-3$\!\pm\!$2.52e-5 & 6.43e-2$\!\pm\!$6.60e-3 & 5.72e-2$\!\pm\!$2.90e-3 & 4.24e-2$\!\pm\!$4.42e-3 & 3.60e-2$\!\pm\!$1.80e-3 \\
\midrule
\multirow{15}{*}{20\%} & \multirow{5}{*}{BIF} & IL1 & 2.79e-3$\!\pm\!$5.33e-5 & \bft{2.27e-3}$\!\pm\!$\bft{9.51e-5} & \bft{7.51e-3}$\!\pm\!$\bft{1.58e-4} & \bft{1.03e-2}$\!\pm\!$\bft{2.40e-4} & \bft{7.42e-3}$\!\pm\!$\bft{7.97e-5} & \bft{9.18e-3}$\!\pm\!$\bft{1.93e-4} \\
 &  & IL5 & \bft{2.76e-3}$\!\pm\!$\bft{2.01e-5} & 2.34e-3$\!\pm\!$1.32e-4 & 7.70e-3$\!\pm\!$2.70e-4 & 1.05e-2$\!\pm\!$2.96e-4 & 7.53e-3$\!\pm\!$1.93e-4 & 9.33e-3$\!\pm\!$2.81e-4 \\
 &  & IL9 & 2.81e-3$\!\pm\!$6.66e-5 & 2.42e-3$\!\pm\!$1.05e-4 & 7.98e-3$\!\pm\!$3.23e-4 & 1.07e-2$\!\pm\!$2.90e-4 & 7.75e-3$\!\pm\!$2.12e-4 & 9.57e-3$\!\pm\!$2.65e-4 \\
 &  & PS & 3.20e-3$\!\pm\!$3.23e-4 & 3.03e-3$\!\pm\!$4.68e-4 & 7.82e-3$\!\pm\!$2.50e-4 & 1.06e-2$\!\pm\!$2.48e-4 & 8.18e-3$\!\pm\!$6.49e-4 & 9.71e-3$\!\pm\!$2.85e-4 \\
 &  & SF & 3.48e-3$\!\pm\!$1.19e-4 & 2.45e-3$\!\pm\!$7.67e-5 & 7.79e-2$\!\pm\!$8.30e-3 & 7.70e-2$\!\pm\!$1.06e-2 & 5.15e-2$\!\pm\!$5.63e-3 & 4.71e-2$\!\pm\!$5.64e-3 \\
\cmidrule{2-9}
 & \multirow{5}{*}{END} & IL1 & \bft{1.10e-3}$\!\pm\!$\bft{1.96e-5} & 4.24e-3$\!\pm\!$3.41e-4 & \bft{5.27e-2}$\!\pm\!$\bft{7.06e-4} & \bft{6.31e-2}$\!\pm\!$\bft{8.41e-4} & \bft{3.72e-2}$\!\pm\!$\bft{5.23e-4} & \bft{4.13e-2}$\!\pm\!$\bft{6.30e-4} \\
 &  & IL5 & 1.15e-3$\!\pm\!$5.99e-5 & 4.00e-3$\!\pm\!$2.34e-4 & 5.48e-2$\!\pm\!$2.15e-3 & 6.51e-2$\!\pm\!$1.35e-3 & 3.86e-2$\!\pm\!$1.49e-3 & 4.27e-2$\!\pm\!$1.08e-3 \\
 &  & IL9 & 1.18e-3$\!\pm\!$1.90e-4 & 4.02e-3$\!\pm\!$2.59e-4 & 5.74e-2$\!\pm\!$9.23e-3 & 6.72e-2$\!\pm\!$7.88e-3 & 4.01e-2$\!\pm\!$5.80e-3 & 4.39e-2$\!\pm\!$5.26e-3 \\
 &  & PS & 1.14e-3$\!\pm\!$2.45e-5 & 4.71e-3$\!\pm\!$5.81e-4 & 5.27e-2$\!\pm\!$7.12e-4 & 6.33e-2$\!\pm\!$8.57e-4 & 3.72e-2$\!\pm\!$5.46e-4 & 4.15e-2$\!\pm\!$6.88e-4 \\
 &  & SF & 2.69e-3$\!\pm\!$1.02e-4 & \bft{3.74e-3}$\!\pm\!$\bft{4.90e-5} & 7.48e-2$\!\pm\!$9.74e-3 & 8.06e-2$\!\pm\!$3.14e-3 & 5.21e-2$\!\pm\!$5.38e-3 & 5.20e-2$\!\pm\!$1.70e-3 \\
\cmidrule{2-9}
 & \multirow{5}{*}{PMU} & IL1 & 2.22e-3$\!\pm\!$6.78e-5 & 3.02e-3$\!\pm\!$6.27e-5 & 5.32e-2$\!\pm\!$3.06e-4 & 6.38e-2$\!\pm\!$4.03e-4 & 3.77e-2$\!\pm\!$2.15e-4 & 4.13e-2$\!\pm\!$2.68e-4 \\
 &  & IL5 & 2.17e-3$\!\pm\!$4.14e-5 & 2.91e-3$\!\pm\!$5.78e-5 & \bft{5.31e-2}$\!\pm\!$\bft{2.69e-4} & \bft{6.37e-2}$\!\pm\!$\bft{3.65e-4} & \bft{3.77e-2}$\!\pm\!$\bft{1.90e-4} & \bft{4.13e-2}$\!\pm\!$\bft{2.35e-4} \\
 &  & IL9 & \bft{2.14e-3}$\!\pm\!$\bft{1.82e-5} & \bft{2.88e-3}$\!\pm\!$\bft{9.08e-6} & 5.33e-2$\!\pm\!$2.87e-4 & 6.39e-2$\!\pm\!$3.68e-4 & 3.78e-2$\!\pm\!$2.07e-4 & 4.14e-2$\!\pm\!$2.46e-4 \\
 &  & PS & 2.42e-3$\!\pm\!$1.89e-4 & 3.26e-3$\!\pm\!$1.12e-4 & 5.32e-2$\!\pm\!$3.65e-4 & 6.38e-2$\!\pm\!$4.54e-4 & 3.79e-2$\!\pm\!$2.52e-4 & 4.14e-2$\!\pm\!$2.99e-4 \\
 &  & SF & 3.02e-3$\!\pm\!$1.13e-4 & 3.75e-3$\!\pm\!$6.23e-5 & 8.12e-2$\!\pm\!$6.06e-3 & 8.25e-2$\!\pm\!$1.98e-3 & 5.51e-2$\!\pm\!$4.02e-3 & 5.31e-2$\!\pm\!$1.34e-3 \\
\bottomrule
\end{tabular}
}
\end{table*}

\bibliographystyle{IEEEtran}
\bibliography{references}


\end{document}